\documentclass[a4paper, 10pt]{amsart}

\usepackage{amsmath,amsthm,amssymb,graphicx,pinlabel,tikz,listings,xcolor}
\usepackage[all]{xy}
\usepackage{tikz-cd}
\usepackage{mathtools}
\usepackage{extarrows}
\usepackage{caption}
\usepackage{subcaption}
\definecolor{wqwqwq}{rgb}{0,0,0}

\SelectTips{cm}{}

\allowdisplaybreaks

\usepackage{hyperref}
\hypersetup{colorlinks=true,citecolor=specialblue,linkcolor=specialblue,urlcolor=black,filecolor=black}

\theoremstyle{definition}
\newtheorem{definition}{Definition}[section]
\newtheorem{remark}[definition]{Remark}

\newtheorem{claim}[definition]{Claim}
\newtheorem{notation}[definition]{Notation}
\newtheorem{question}[definition]{Question}
\newtheorem{example}[definition]{Example}

\theoremstyle{plain}
\newtheorem{theorem}[definition]{Theorem}
\newtheorem{proposition}[definition]{Proposition}
\newtheorem{lemma}[definition]{Lemma}
\newtheorem{corollary}[definition]{Corollary}

\makeatletter
    
    \@addtoreset{equation}{section}%
  \makeatother

\makeatletter
\def\blfootnote{\xdef\@thefnmark{}\@footnotetext}

\definecolor{specialblue}{RGB}{35,88,180}
\definecolor{lightblue}{RGB}{10,180,255}
\definecolor{verylightblue}{RGB}{213,255,255}

\newcommand{\Q}{\mathbb{Q}}
\newcommand{\Z}{\mathbb{Z}}

\newcommand{\calI}{\mathcal{I}}
\newcommand{\calIC}{\mathcal{IC}}
\newcommand{\calM}{\mathcal{M}}
\newcommand{\calT}{\mathcal{T}}
\newcommand{\calS}{\mathcal{S}}

\newcommand{\calG}{\mathcal{G}}
\newcommand{\calC}{\mathcal{C}}

\newcommand{\calL}{\mathcal{L}}
\newcommand{\calN}{\mathcal{N}}

\newcommand{\Aut}{\operatorname{Aut}}
\newcommand{\IAut}{\operatorname{IAut}}

\newcommand{\Tr}{\operatorname{Tr}}
\newcommand{\Ker}{\operatorname{Ker}}
\newcommand{\loga}{\operatorname{log}}
\newcommand{\mir}{\operatorname{Mir}}
\newcommand{\Gr}{\operatorname{Gr}}

\newcommand{\ev}{\operatorname{ev}}

\newcommand{\im}{\operatorname{Im}}
\newcommand{\Der}{\operatorname{Der}}

\newcommand{\Sp}{\operatorname{Sp}}

\newcommand{\GL}{\operatorname{GL}}

\newcommand{\trait}{\!\!\begin{array}{c}{\rule{5mm}{0.4mm}}\\[0.1cm] \end{array}\!\!}

\newcommand{\vltree}[4]{\tikz[baseline=6pt, scale = 0.4]{
\draw [color=wqwqwq] (0,2-1.25)-- (2,2-1.25);
\draw [color=wqwqwq] (0,1.5)-- (0,0);
\draw [color=wqwqwq] (2,1.5)-- (2,0);

\draw[color=wqwqwq] (0,3-1.15) node {$#1$};
\draw[color=wqwqwq] (0,1.1-1.5) node {$#2$};
\draw[color=wqwqwq] (2,1.1-1.5) node {$#3$};
\draw[color=wqwqwq] (2,3-1.15) node {$#4$};}}

\newcommand{\vlfivetree}[5]{\tikz[baseline=6pt, scale = 0.4]{
\draw [color=wqwqwq] (0,2-1.25)-- (2.5,2-1.25);
\draw [color=wqwqwq] (0,1.5)--(0,0);
\draw [color=wqwqwq] (2.5,1.5)--(2.5,0);
\draw [color=wqwqwq] (1.25,0.75)--(1.25,1.5);

\draw[color=wqwqwq] (0,3-1.1) node {$#1$};
\draw[color=wqwqwq] (0,1.1-1.5) node {$#2$};
\draw[color=wqwqwq] (2.5,1.1-1.5) node {$#3$};
\draw[color=wqwqwq] (2.5,3-1.1) node {$#4$};
\draw[color=wqwqwq] (1.25,3-1.1) node {$#5$};}}

\newcommand{\ltritree}[3]{\tikz[baseline=8pt ,scale = 0.5]{
\draw [color=black] (1,1.5)-- (1,0.75);
\draw [color=black] (1,0.75)-- (1-0.866*0.75,0.75-.37);
\draw [color=black] (1,0.75)-- (1+0.866*0.75,0.75-.37);
\draw[color=black] (1,1.5+0.3) node {$#1$};
\draw[color=black] (1-0.866*0.75-0.2,0.75-.37-0.2) node {$#2$};
\draw[color=black] (1+0.866*0.75+0.2,0.75-.37-0.2) node {$#3$};}}

\newcommand{\rtritree}[3]{\tikz[baseline=10pt, scale = 0.5]{
\draw[color=wqwqwq] (0.75,0)-- (0.75,0.75);
\draw[color=wqwqwq] (0.75,0.75) -- (0,1.5);
\draw[color=wqwqwq] (0.75,0.75) --(1.5,1.5);
\node[label= {[xshift=0cm, yshift=0.6cm]{$#1$}}]{} ;
\node[label= {[xshift=0.75cm, yshift=0.6cm]{$#2$}}]{} ;
\node[label= {[xshift=0.375cm, yshift=-0.6cm]{$#3$}}]{} ;} }

\newcommand{\treethreel}[3]{\tikz[baseline=10pt, scale = 0.8]{
\draw[color=wqwqwq] (0.5,0)--(0.5,0.25);
\draw[color=wqwqwq] (0.5,0.25)--(0,0.75);
\draw[color=wqwqwq] (0.25,0.5)--(0.5,0.75);
\draw[color=wqwqwq] (0.5,0.25)--(1,0.75);
\draw[color=wqwqwq] (0,1) node {\small $#1$};
\draw[color=wqwqwq] (0.5,1) node {\small $#2$};
\draw[color=wqwqwq] (1,1) node {\small $#3$};}}

\newcommand{\treethreer}[3]{\tikz[baseline=10pt, scale = 0.8]{
\draw (0.5,0)--(0.5,0.25);
\draw (0.5,0.25)--(0,0.75);
\draw (0.75,0.5)--(0.5,0.75);
\draw (0.5,0.25)--(1,0.75);
\draw[color=wqwqwq] (0,1) node {\small $#1$};
\draw[color=wqwqwq] (0.5,1) node {\small $#2$};
\draw[color=wqwqwq] (1,1) node {\small $#3$};}}

\newcommand{\btreethreel}[3]{\tikz[baseline=10pt, scale = 1.2]{
\draw[color=wqwqwq] (0.5,0)--(0.5,0.25);
\draw[color=wqwqwq] (0.5,0.25)--(0,0.75);
\draw[color=wqwqwq] (0.25,0.5)--(0.5,0.75);
\draw[color=wqwqwq] (0.5,0.25)--(1,0.75);
\draw[color=wqwqwq] (0,1) node {\small $#1$};
\draw[color=wqwqwq] (0.5,1) node {\small $#2$};
\draw[color=wqwqwq] (1,1) node {\small $#3$};}}

\newcommand{\btreethreer}[3]{\tikz[baseline=10pt, scale = 1.2]{
\draw (0.5,0)--(0.5,0.25);
\draw (0.5,0.25)--(0,0.75);
\draw (0.75,0.5)--(0.5,0.75);
\draw (0.5,0.25)--(1,0.75);
\draw[color=wqwqwq] (0,1) node {\small $#1$};
\draw[color=wqwqwq] (0.5,1) node {\small $#2$};
\draw[color=wqwqwq] (1,1) node {\small $#3$};}}

\newcommand\hooklongrightarrow{\mathrel{\lhook\mkern -3.5mu\relbar\mkern -4.5mu \rightarrow }}

\author{Quentin Faes}
\newcommand{\adresse}{{
  \bigskip
  \footnotesize
  \textsc{Institut fur Mathematik, Universitat Zurich, CH-8057 Zurich}\par\nopagebreak
  \textit{E-mail address}  \texttt{quentin.faes@math.uzh.ch}}}

\begin{document}
\begin{abstract}
The so-called Johnson homomorphisms $(\tau_k)_{k \geq 1}$ embed the graded space associated to the Johnson filtration of the mapping class group of a surface with one boundary component into the Lie ring of positive symplectic derivations $D(H)$. In this paper, we show the existence of torsion in the cokernels of the Johnson homomorphisms, for all even degrees, provided the genus is big enough. For each degree, we define a map on a subgroup of $D(H)$, whose image is 2-torsion, and which vanishes on the image of $\tau$. We also give a formula to compute these maps.
\end{abstract}

\title{Torsion in the cokernels of the Johnson homomorphisms}
\maketitle

\vspace{-0.5cm}
 \tableofcontents

\blfootnote{This research was funded by the Japanese Society for the Promotion of Science via a JSPS Postdoctoral Fellowship for Research in Japan (Graduate school of Mathematical Sciences, University of Tokyo).}

\section{Introduction}

The mapping class group of a surface with one boundary component can be equipped with the so-called Johnson filtration, whose associated graded space is a Lie ring that can be embedded in a Lie ring of derivations via the \emph{Johnson homomorphisms}. While most results concerning these homomorphisms are known over the rationals, we shall focus here on some torsion questions, which require a more subtle approach.

We consider a compact connected oriented surface $\Sigma_{g,1}$ of genus $g$ with one boundary component, and a point $\star$ in $\partial \Sigma$. We denote by $\pi := \pi_1(\Sigma_{g,1},\star)$ its fundamental group and by $H : = H_1(\Sigma_{g,1}, \Z)$ its first homology group. Setting $n = 2g$, the group $\pi$ is isomorphic to $F_n$, the free group on $n$ generators, and $H$ is a free abelian group on $n$ generators. For an element $x \in \pi$, we shall write $\{ x \}$ for the class of $x$ in $H$. We denote by $\calM := \calM(\Sigma_{g,1})$ the mapping class group of the surface $\Sigma_{g,1}$. The group $\calM$ acts naturally on $\pi$ and we denote this action by $$\rho : \calM \longrightarrow \Aut(\pi).$$

By the Dehn-Nielsen-Baer theorem, the representation $\rho$ is faithful, and its image coincides with the subgroup of automorphisms of $\pi$ that preserve  $\zeta \in \pi$, the class of the oriented boundary component (see for example \cite[Theorem 8.1]{fm}). This motivates the study of $\calM$ as a subgroup of $\Aut(F_n)$.

We define $A_1$ (resp. $J_1$) as the subgroup of $\Aut(F_n)$ (resp. $\calM$) acting trivially on $H$. The group $J_1$ is otherwise known as the \emph{Torelli group} $\calI$ of the surface, i.e. the subgroup of mapping classes acting trivially on homology, and $A_1$ is usually denoted by $\operatorname{IAut}(F_n)$. To be more precise, the action of $\calM$ on $H$ preserves the symplectic form $\omega : H\otimes H \rightarrow \Z$ induced by the intersection form of $\Sigma_{g,1}$. Using the generators of the symplectic group from \cite{klingen} or \cite{mennicke}, one can show that if $\Sp(H)$ denotes the symplectic group of $H$, then $\calM \rightarrow \Sp(H)$ is surjective (see for example \cite{birman75}). In summary, we get the following commutative diagram of exact sequences: 

\begin{equation}
\label{hactions}
\begin{tikzcd} 
0  \arrow[r]  & \calI \arrow[r] \arrow[d, hook] & \calM \arrow[r] \arrow[d, hook] & \Sp(H) \arrow[r] \arrow[d,hook] & 0 \\ 0  \arrow[r] & A_1 \arrow[r] & \Aut(F_n) \arrow[r] & \operatorname{GL}(H)\arrow[r] & 0.
\end{tikzcd}
\end{equation}

For a group $G$, we denote by $(\Gamma_kG)_{k \geq 1}$ the lower central series of $G$ (defined inductively by $\Gamma_1G = G$ and $\Gamma_{k+1}G = [G, \Gamma_kG]$). The $k$-th nilpotent quotient of $G$ is, by definition, $G/ \Gamma_{k+1}G$. By definition $H = F_n / \Gamma_2 F_n$, hence we could have defined $A_1$ as the kernel of the action of $\Aut(F_n)$ on the first nilpotent quotient of $F_n$. Similarly, the \emph{Andreadakis-Johnson filtration} $(A_k)_{k\geq1}$ of $\Aut(F_n)$ is defined by setting $A_k$ to be the kernel of the action of $\Aut(F_n)$ on the $k$-th nilpotent quotient of $F_n$. One can then define the Johnson filtration $(J_k)_{k \geq1}$ of $\calM$ by $J_k := \calM \cap A_k$. These two filtrations are important, as they induce separated topologies on $\Aut(F_n)$ and $\calM$ (meaning that the intersection of all the terms of the filtration is trivial). This is due to the fact that free groups are residually nilpotent.

Hence $\calM$ can be understood by studying the ``slices" $J_k/J_{k+1}$ of the filtration. This is done by using the technology of Johnson homomorphisms. We set $\Gr^A(\Aut(F_n)) := \bigoplus A_k/A_{k+1}$ and $\Gr^J(\calM) := \bigoplus J_k/J_{k+1}$. These graded spaces have a Lie ring structure whose Lie bracket is induced by taking group commutators in $\calM$ and $\Aut(F_n)$. Moreover, there is a well-defined injective graded Lie ring morphism, called the \emph{Johnson homomorphism}:
 \[  \tau: \Gr^J(\mathcal{M}) \longrightarrow D(H) \] where the target space is some Lie ring of derivations. The main result of this paper is the following.\\

\noindent \textbf{Theorem A.} \emph{For any $k \geq 1$, for any $g \geq k+2$, $D_{2k}(H) / \im(\tau_{2k})$ has non-trivial $2$-torsion.} \\

Before reviewing the definitions more precisely, let us mention why Theorem A can be of interest for low-dimensional topologists. First, the Johnson filtration is involved in the definition of some equivalence relations on 3-manifolds that are not yet understood properly, and that are related to finite-type invariants of homology 3-spheres \cite{masY3, faes2}. As an illustration of this fact, the Johnson homomorphisms can be recovered from the tree-level of the LMO invariant \cite{masinf}. Recently, the Johnson filtration also appeared in the descriptions or partial descriptions of the kernels of some quantum representations \cite{DetSan} and some homological representations \cite{mori,BPS}. Finally, the Johnson homomorphisms also extend to the monoid of homology cylinders, and thus can be used to study the difference between the mapping class group and this monoid. This is illustrated in \cite{sakgoda}, where it is claimed that the Johnson homomorphisms can help to detect non-fibered homologically fibered knots. See also \cite{NSSred} where the Johnson homomorphisms (and the Satoh trace defined below) appear in the 1-loop part of the LMO homomorphism.

We now review the definition of $\tau$ by giving a list of properties, that can be found scattered in papers of Johnson in the 80's \cite{joh80,johsurvey} and Morita in the 90's \cite{mor91, mor93}. Denote by $\calL(H)$ the graded free Lie ring generated by $H$ in degree 1, we have $\calL_{k+1}(H) \simeq \Gamma_{k+1}F_n/\Gamma_{k+2}F_n$, for any $k \geq 0$. We set $\Der_k(H)$ the Lie ring of derivations of $\calL(H)$ that shifts the degree by $k$. If $k$ is positive, then the derivation is said to be \emph{positive}. The process of restricting a derivation to $H = \calL_1(H)$ provides an identification $\Der_k(H) \simeq H^* \otimes \calL_{k+1}(H)$. 

\begin{itemize}
\item For any $k,l \geq 1$, we have $[A_k, A_l] \subset A_{k+l}$ and $[J_k, J_l] \subset J_{k+l}$. In particular $[A_k,A_k] \subset A_{2k} \subset A_{k+1}$ and $A_k/ A_{k+1}$ is abelian. Similarly $J_k / J_{k+1}$ is abelian.
\item For any $x\in \pi$, and $f \in A_k$, the class of $f(x)x^{-1}$ in $\Gamma_{k+1}\pi/\Gamma_{k+2}\pi$ depends only on the homology class of $x$. Hence we get a well defined map $\tau_k' : A_k \rightarrow \Der_k(H)$, by the formula \[ \tau'_k(f)(\{ x\}) = [f(x)x^{-1}].\]
\item This yields a short exact sequence of groups \[ 
\begin{tikzcd}
0 \arrow[r] & A_{k+1} \arrow[r] & A_k \arrow[r, "\tau_k'"] & \Der_k(H). 
\end{tikzcd}\]
\item $\Aut(F_n)$ acts on itself by conjugation, inducing a $\GL(H)$-action on $A_k / A_{k+1}$, and $\tau_k' : A_k/ A_{k+1} \rightarrow \Der_k(H)$ is $\GL(H)$-equivariant (here, the action of $f_* \in \GL(H)$ on a derivation $d$ is given by $f_* \cdot d := f_* \circ d \circ f_*^{-1}$). 
\end{itemize}

We shall call $\tau_k'$ the \emph{$k$-th Andreadakis-Johnson homomorphism}. It is known (see \cite{mor93}) that for an element in $J_k$, $\tau_k'(f)$ sends $\omega$, the element of $ \calL_2(H)$ induced by $\zeta^{-1}$, to $0$ (see also Section \ref{sec3}). We denote by $D_k(H)$ the subset of derivations that send $\omega$ to $0$. Such a derivation is called a \emph{symplectic} derivation. Hence, we can define, by restriction and corestriction of $\tau_k'$, the $k$-th Johnson homomorhism $\tau_k$:
\begin{align*}
\tau_k: J_k & \longrightarrow D_k(H) \\
f & \longmapsto \big( \{ x\} \mapsto [f(x)x^{-1}] \in \calL_{k+1}(H) \big).
\end{align*}

\noindent $\calM$ acts on itself by conjugation, inducing an $\Sp(H)$-action on the quotients $J_k / J_{k+1}$, making $\tau_k : J_k / J_{k+1} \rightarrow D_k(H)$ a map of $\Sp(H)$-module.

Let $\Gr^A(\Aut(F_n)) := \bigoplus A_k/A_{k+1}$ and $\Gr^J(\calM) := \bigoplus J_k/J_{k+1}$. These graded spaces have a Lie ring structure. According to the discussion above, we have the following commutative diagram of graded Lie ring morphisms: 
\[
\begin{tikzcd}  
\Gr^A(\Aut(F_n)) \arrow[r, "\tau'" ,hook] & \Der(H)  \\
\Gr^J(\mathcal{M}) \arrow[r, "\tau", hook]\arrow[u, hook] & D(H).\arrow[u,hook]
\end{tikzcd}
\]

As explained above, it is an important question to determine the images of $\tau'$ and $\tau$. A derivation is called \emph{tame}, if it is in the image of $\tau'$. Morita \cite{mor93} gave a first obstruction for a symplectic derivation to be in the image of $\tau$. As explained in Section \ref{sec2}, Satoh generalized Morita's work by introducing a trace-like operator $\Tr_k$ on $\Der_k(H)$ whose kernel, for $k \geq 2$ and $n \geq k+2$, is equal to $\im(\tau_k')$ \cite{sat12, darne}. Unfortunately, this result does not hold for $\im(\tau)$, i.e. there are some symplectic tame derivations which are not in the image of $\tau$. For example, Satoh's trace vanishes on $D_2(H)$, but the quotient $D_2(H)/ \im(\tau_2)$ is non-trivial. The image of $\tau_2$ is described by Morita and Yokomizo in \cite{mor} and \cite{yok}. The cokernel $D_2(H) / \im(\tau_2)$ is a non-trivial torsion group, and the $\Sp(H)$-module structure of the cokernel is described explicity by the author in \cite{faes}.

The restriction of the Satoh Trace to symplectic derivations was studied by Enomoto and Satoh in \cite{ES}. Then, Conant \cite{conant}, building on \cite{CKV}, defined a trace $\Tr^C$ which generalizes $\Tr$ and vanishes on $\im(\tau)$. But this result is based on a theorem of Hain \cite{haininf} implying that, rationally, the image of $\tau_k$ is generated in degree 1. Thus, except for degree 2, no information was known, as far as the author knows, about torsion in the cokernels of the Johnson homomorphisms. In this paper we define a new trace map $\overline{\Tr}$ defined on $\im(\tau'_k) \cap D_k(H)$ which detects torsion in every even degree.\\

We set $T(H)$ the tensor algebra generated by $H$ in degree $1$ and $C(H):= T(H) / [T(H), T(H)]$ the space of cyclic tensors (here for $x,y \in T(H)$, $[x,y]:= xy - yx$). The Satoh trace is valued in $C(H)$. In Section \ref{sec3}, we set $\overline{C(H)} := \Tr(D_k(H))$ and in Section \ref{sec6} we define $\mir_k(H)$ to be the subspace of $C_k(H)$ generated by elements of the form $w + (-1)^{k+1}\overline{w}$, for $w \in T_k(H)$ where $\overline{w}$ is the mirror image of the tensor $w$. Theorem A will follow from Theorem B.\\

\noindent \textbf{Theorem B.} \emph{Let $k\geq 2$ and $g \geq k+2$. There exists a well-defined graded $\Sp(H)$-linear map $$\overline{\Tr}_k : \im(\tau'_k) \cap D_k(H) \rightarrow C_{k+1}(H)/\overline{C_{k+1}(H)}$$ vanishing on $\im(\tau_k)$, non-trivial for $k$ even, and whose image consists of 2-torsion elements.} \\

\noindent For the proof of the non-triviality in Theorem B, we shall prove that there exists a map $$\Tr^{\mir}_k :  \im(\tau'_k) \cap D_k(H) \rightarrow C_{k+1}(H)/\mir_{k+1}(H), $$ obtained as a quotient of $\overline{\Tr}$, which is non-trivial in even degrees. This map vanishes when $k$ is odd, though.\\

The heuristic with which the map $\overline{\Tr}$ is the following:\\

\noindent \textbf{Heuristic C.} \label{heuri}
Considering a derivation $d$ of degree $k$, we want to determine if it sits in $\im(\tau_k)$ or not. We have a series of obstructions.

\begin{enumerate}
\item We need to check first that $d \in \Ker(\Tr_k)$, so that $d$ belongs at least stably in the image of $\tau'$. (\textbf{tameness})
\item Then $d(\omega)$ should be $0$, which is equivalent (see Claim \ref{claim0}) to the fact that any lift $f$ of $d$ to $A_k$ verifies $f(\zeta)\zeta^{-1} \equiv 1 \in \Gamma_{k+2}\pi / \Gamma_{k+3} \pi$, i.e. $d \in D_k(H)$. (\textbf{symplecticity})
\item This process can be carried on, i.e. we can check, among the lifts of $d$ to $A_k$, if there exists one such that $f(\zeta)\zeta^{-1} \equiv 1 \in \Gamma_{k+3}\pi / \Gamma_{k+4}\pi$.
\item If this is the case, then we can look \emph{among the lifts $f$ verifying \text{(3)}} if there is one such that $f(\zeta)\zeta^{-1} \equiv 1 \in \Gamma_{k+4}\pi / \Gamma_{k+5}\pi$, and so forth.
\end{enumerate}

\noindent In this paper we shall focus on (3) and we shall obtain the map $\overline{\Tr}$. Theoretically, we can define ``higher" obstructions. Nevertheless, it is not clear how to compute those. We discuss this in Remark \ref{finalremark}, as well as the efficacy of these obstructions.

The paper is organized as follows. In Section \ref{sec2}, we recall the definitions of the Andreadakis and Johnson filtrations as well as the Satoh Trace. Section \ref{sec3} is dedicated to the introduction of the obstruction map $\overline{\Tr}$, following (3) in Heuristic C. In Section \ref{sec4}, we introduce the so-called infinitesimal Andreadakis-Johnson representation and relate it to the compuation of the map $\overline{\Tr}$. In Section \ref{sec5}, we use this relation and some properties of the infinitesimal representation to study the map $\overline{\Tr}$. In Section \ref{sec6}, we study the target space of $\overline{\Tr}$ and introduce the map $\Tr^{\mir}$. We also provide some examples, and relate our work to results previously known about the cokernel of $\tau_2$. In Section \ref{sec7}, we provide a formula to compute $\overline{\Tr}$, and use it to give explicit examples of torsion elements in the cokernel of $\tau_k$, for $k$ even. In the last section, we give a topological interpretation for $\overline{\Tr}$ in degree $2$, and briefly discuss potential generalizations.

\vspace{4mm}
\noindent
\textbf{Acknowledgment.} The author would like to thank T. Sakasai for hosting him in the University of Tokyo, and G. Massuyeau for suggesting the computations in Section \ref{sec6}. The author is also thankful to both of them for making comments on the earliest version of this paper, and to an anonymous referee for suggesting the current reorganization of this paper.

\section{About the Andreadakis and Johnson filtrations}
\label{sec2}
We first review the definition and properties of the Satoh trace $\Tr$, following \cite{sat12}, \cite{darne} and \cite{massak}. We then recall the definition of the Johnson filtration of the mapping class group and some facts about the restriction of the Andreadakis-Johnson homomorphisms to this filtration. 

\subsection{Satoh's trace-like operator}

Let $T(H)$ be the the tensor algebra generated by $H$ in degree 1, and $[T(H),T(H)]$ the subspace spanned by commutators $[x,y] := xy -yx$ for $x,y \in T(H)$. Let $C(H) := T(H) / [T(H), T(H)]$ be the space of cyclic tensors. This quotient inherits the graduation of $T(H)$. It also inherits an action of $\GL(H)$. It is a well-known fact that $\Der(H):= \Der(\calL(H),\calL(H))$, the space of derivations of the the free Lie algebra, is isomorphic as a vector space to $\operatorname{Hom}(H,\calL(H)) = H^* \otimes \calL(H) $. Satoh considered the $\GL (H)$-equivariant graded linear map \[ \Tr : \Der(H) \longrightarrow C(H) \] which is defined in degree $k \geq 0$ by the composition \[ \Der_k(H) = H^* \otimes \calL_{k+1}(H) \hooklongrightarrow H^* \otimes T_{k+1}(H) \longrightarrow T_{k}(H) \longrightarrow C_k(H)  \] where the free Lie algebra $\calL(H)$ naturally imbeds in $T(H)$ and the second arrow is induced by the evaluation of the factor $H^*$ on the first factor in $T_{k+1}(H)$.

In \cite[Theorem 3.1]{sat12}, Satoh used the map $\Tr$ to compute, rationally, the cokernel of the Andreadakis-Johnson homomorphism $\tau_k'$ when $n \geq k+2$. In \cite[Prop. 2.37]{darne}, Darn\'e highlighted that this results holds over the integers.

\begin{theorem}[Satoh, Darn\'e]
\label{satoh}
Let $k \geq 2$ and $n \geq k+2$ be integers, and $\tau_k': A_k(F_n) \rightarrow \Der_k(H)$ be the $k$-th Andreadakis-Johnson homomorphism, then \[\im( \tau'_k ) = \Ker (\Tr_k)\] and \[\operatorname{Coker}(\tau_k')= C_k(H).\]
\end{theorem}  

\noindent In other words, a derivation is tame, up to stabilization, if and only if it is in the kernel of the Satoh trace. By stabilizing, we mean increasing the rank $n$.

We now review some results from \cite{massak} that will be useful in the sequel. Note that any derivation of $\calL (H)$ extends uniquely to a derivation of $T(H) $, preserving $[T(H), T(H)]$ (since $d([x,y]) = d(xy - yx) = [d(x),y] + [x,d(y)]$). Thus, such derivations induce linear endomorphisms of $C(H)$, defining a $\Der(H)$-module structure on $C(H)$. The map $\Tr$ is a 1-cocycle on the Lie algebra $\Der(H)$ with respect to this action: 

\begin{proposition}[Massuyeau-Sakasai]
\label{cocyle}
For any $\delta, \eta \in \Der(H)$, we have \[ \Tr([\delta, \eta]) = \delta \cdot \Tr(\eta) - \eta \cdot \Tr(\delta). \]
\end{proposition}

\begin{remark}
\label{remarkcocycle}
By definition, the action of a positive derivation on an element in $C_1(H)$ is trivial. Proposition \ref{cocyle} then implies among other things that if $\delta$ is a positive derivation such that $\Tr(\delta) = 0$ and $\eta$ is of degree 1, then $\Tr([\delta, \eta]) = 0$. We will make use of this fact later. This also proves that, in degree strictly greater than $1$, the subalgebra of $\Der(H)$ generated in degree $1$ is in the kernel of the trace.

\end{remark}
A fundamental property is the following result \cite[Prop. 5.3]{massak}. We will refine it below using the proof given by its authors (Prop \ref{massakimproved}).

\begin{proposition}[Massuyeau-Sakasai]
\label{massak}
For any $n \geq 2$ and $k \geq 2$, the map $\Tr_k$ vanishes on the image of $\tau_k'$.
\end{proposition}

\subsection{The Johnson filtration and the Johnson homomorphisms}
\begin{figure}[h]
\includegraphics[scale=0.8]{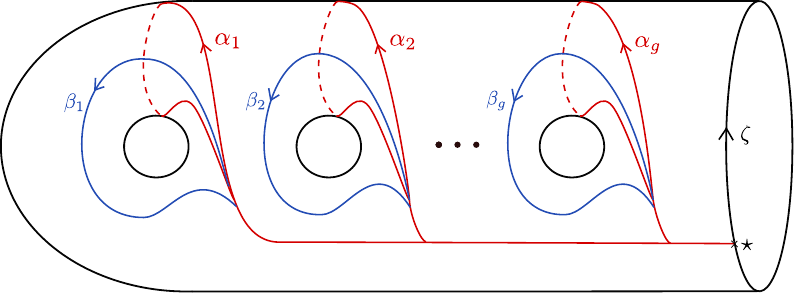}
\caption{A convenient choice of basis for $\pi$}
\label{figsurf}
\end{figure}

In this section, we consider a surface $\Sigma_{g,1}$ of genus $g$ with one boundary component. The fundamental group $\pi$ of the surface is identified with  $F_n$ for $n = 2g$, by choosing a basis of the fundamental group denoted by $\alpha_1, \beta_1, \alpha_2, \beta_2, \cdots, \alpha_g, \beta_g$. The classes of $\alpha_i$ and $\beta_i$ in $H$ (identified with the first homology group of the surface) are denoted by $a_i$ and $b_i$ and the basis is chosen so that $\omega(a_i,b_i) = 1$. Recall that $\omega$ is the symplectic form given by the intersection form on the surface. Figure \ref{figsurf} gives an example of such a choice of basis. 

The mapping class group $\calM$ of the surface acts faithfully on $\pi$ by the Dehn-Nielsen theorem. Via this representation $\rho$, the group $\mathcal{M}$ coincides with the subgroup of automorphisms that fixes the element $\zeta \in \Gamma_2\pi$ corresponding to the boundary component. We also denote by $\omega := \sum_{i=1}^{g} [a_i,b_i]$ the class of $\zeta^{-1} \in \calL_2(H)$. For example, with the specific basis chosen in Figure \ref{figsurf}, we have $\zeta^{-1} := \prod_{i=1}^g [\beta_i^{-1}, \alpha_i]$. The Andreadakis filtration restricts to the Johnson filtration: \[ J_k = A_k \cap \calM\] and the Andreadakis-Johnson homomorphism $\tau_k'$ restricts to the Johnson homomorphism \[ \tau_k : J_k \longrightarrow \Der_k(H)\] which is equivariant with respect to the mapping $\calM \rightarrow \Sp(H)$ induced by the action of $\mathcal{M}$ on $H$. The following is well-known, but is also proven below in Prop. \ref{importantlemma}.

\begin{claim}
\label{claim0}
It is known that for any $f \in A_k$, we have $f(\zeta)\zeta^{-1} \in \Gamma_{k+2}\pi$, and \[ f(\zeta)\zeta^{-1} \equiv -\tau_k(f)(\omega) \in \calL_{k+2}(H). \]
\end{claim}

Hence for any element $f$ in $J_k$, $\tau_k(f)(\omega)$ is zero. Let $\ev_\omega : \Der_k(H) \rightarrow \calL_{k+2}(H)$ be the evaluation at $\omega$ and let $D(H)$ be defined by the following exact sequence:

\begin{equation}
\label{ses}
0 \longrightarrow D_k(H) \longrightarrow \Der_k(H) \xlongrightarrow{\ev_\omega} \calL_{k+2}(H) \longrightarrow 0.  
\end{equation}

\noindent $D(H)$ is the so-called space of \emph{symplectic derivations}. For further use, we prove the surjectivity of $\ev_\omega$.

\begin{proposition}
\label{propses}
The short sequence \eqref{ses} is exact.
\end{proposition}

\begin{proof}

It suffices to find lifts for elements of the form $[a_i,u_i]$ and $[b_i,v_i]$ for $u_i, v_i \in \calL_{k+1}(H)$ brackets of elements of the basis of $H$ given by the $a_i$'s and the $b_i's$. Indeed, these brackets linearly generate $\calL_{k+2}(H)$. For a derivation $d$, we have $d(\omega) = \sum [d(a_i),b_i] + [a_i,d(b_i)]$. Hence we can pick  $b_i^* \otimes u_i$ for $[a_i,u_i]$ and $-a_i^* \otimes v_i$ for $[b_i,v_i]$.
\end{proof} 

By definition of $D_k(H)$, $\tau_k$ maps $J_k$ to $D_k(H)$, hence the following commutative diagram follows: 

\[
\begin{tikzcd}  
A_k \arrow[r, "\tau_k'"] & \Der_k(H)  \\
J_k \arrow[r, "\tau_k"]\arrow[u, hook] & D_k(H).\arrow[u,hook]
\end{tikzcd}
\]

From now on, we consider the map $\tau_k$ to be valued in $D_k(H)$ and the goal of this paper is to study its cokernel. Sometimes, in our notation, $\tau_k$ also refers to the induced homomorphism of abelian group $J_k /J_{k+1} \rightarrow D_k(H)$ and we set $\tau_k^\Q$ to be the rationalization of this map.

\begin{remark}
Some tame symplectic derivations are not in the image of $\tau_k$, even rationally. Indeed, Conant, Kassabov and Vogtmann \cite{CKV} and Conant \cite{conant} showed that the image of $\tau_k^\Q$ is smaller than the intersection of the kernel of $\Tr^\Q$ and $D_k(H)$.
\end{remark}

\subsection{Cokernel of $\tau_2$}

We now recall what is known in the case $k=2$. We will later generalize to define new maps on the cokernels of the Johnson homomorphisms. In the sequel, $a_i^*$ (resp $b_i^*$)  denotes the linear map sending $a_i$ (resp. $b_i$) to $1$, and the other elements to $0$. It is important to note that $a_i^*$ is not the dual of $a_i$ with respect to the symplectic form $\omega$ (this dual would be $b_i^* = \omega(a_i,-)$).

The cokernel of $\tau_2$ was computed in \cite{mor} and \cite{yok}. Morita proved it is 2-torsion, and Yokomizo computed the rank over $\mathbb{Z}_2$. In \cite{faes}, the author introduced a trace-like operator $\Tr^{as} : D_2(H) \rightarrow \Lambda^2(H) \otimes \mathbb{Z}_2$ such that the following short sequence is exact.

\begin{equation}
\begin{tikzcd}
0 \arrow[r]& \im{(\tau_2)} \arrow[r]& D_2(H) \arrow[r, "\Tr^{as}"]& \Lambda^2(H) \otimes \mathbb{Z}_2 \arrow[r] &0.
\end{tikzcd}
\end{equation}

\begin{example}
\label{exampled}
The derivation $d := b_1^* \otimes [a_2,[a_1,a_2]] + b_2^* \otimes [a_1,[a_2,a_1]]$ is symplectic and tame, i.e. $d \in \Ker(\Tr) \cap D_2(H)$. Nevertheless $\Tr^{as}(d) = a_1 \wedge a_2$, hence this derivation is not in the image of $\tau_2$. 
\end{example}

\begin{remark}
The Satoh Trace actually vanishes on $D_2(H)$, hence all symplectic derivations of degree $2$ are tame (for $g \geq 2$).
\end{remark}

\noindent The computations of Morita, Yokomizo, and the author are consequences of the fact that we know generators for the so-called \emph{Johnson kernel} $J_2$ \cite{joh2}, but there was no direct justification of the fact that $\Tr^{as}$ vanishes on $\im(\tau_2)$ (a topological justification of this fact will be given in Section \ref{sec6}, though). We actually claim the following, which is related to Heuristic C in the introduction.

\begin{claim}
\label{claim}
For any $d \in D_2(H)$ such that $\Tr^{as}(d) \neq 0$, and for any lift $f$ of $d$ to $A_2$, we have $f(\zeta)\zeta^{-1} \in \Gamma_{5}\pi$ and its class in $\calL_{5}(H)$ is non-trivial.
\end{claim}

\noindent Of course, this implies that such a derivation is not in the image of $\tau_2$, as otherwise one would have the existence of a lift $f$ with $f(\zeta)\zeta^{-1} = 1$. Claim \ref{claim} will be a consequence of Proposition \ref{tracesrel} and holds, for example, for the derivation $d$ defined in Example  \ref{exampled}.

\section{Definition of the new obstruction}
\label{sec3}

Let $k \geq 1$ and $\calG_k := \{ f \in A_k | f(\zeta) \zeta^{-1} \in \Gamma_{k+3}\pi \} = A_k \cap \tau_k'^{-1}(D_k(H))$ where the equality follows from Claim \ref{claim0}. This is clearly a subgroup of $A_k$ containing $A_{k+1}$ and $J_k$. It is also stable by conjugation by elements of $\calM$.

\begin{proposition}\label{propmumorphism}
The map 
\begin{align*}
\mu_k: \calG_k & \longrightarrow  \mathcal{L}_{k+3}(H) \\
f & \longmapsto f(\zeta) \zeta^{-1} \in \Gamma_{k+3}\pi / \Gamma_{k+4}\pi \cong \mathcal{L}_{k+3}(H).
\end{align*} is a group homomorphism.
\end{proposition}

\begin{proof}
For $f, g \in \mathcal{G}_k$, we have $(f \circ g)(\zeta) \zeta^{-1}=f \left(g(\zeta) \zeta^{-1}\right) f(\zeta) \zeta^{-1}$. But since $f$ belongs to $A_1$ it acts trivially on $\Gamma_{k+3}\pi / \Gamma_{k+4}\pi$, and we have $f\left(g(\zeta) \zeta^{-1}\right) \equiv$ $g(\zeta) \zeta^{-1}$ modulo $\Gamma_{k+4}\pi$. .
\end{proof}

\begin{proposition}\label{propeq1}
The map $\mu_k$ is $\calM$-equivariant for the action of $\calM$ by conjugation on $\calG_k$ and the symplectic action on $\calL_{k+3}(H)$.
\end{proposition}

\begin{proof}
If $g \in \mathcal{M}$, then $g(\zeta)=\zeta$ and $\left(gf g^{-1}\right)(\zeta) \zeta^{-1}=$ $g f(\zeta) \zeta^{-1}=g f(\zeta) g(\zeta)^{-1}=g\left(f(\zeta) \zeta^{-1}\right).$ 
\end{proof}
 
The restriction of $\mu_k$ to $A_{k+1}$ is clearly $f \mapsto \tau_{k+1}'(f)(\omega)$, \emph{i.e.} $\ev_\omega \circ \tau_{k+1}'(f)$. Given $f \in \mathcal{G}_k$, we want to know if we can multiply it by an element $h$ of $A_{k+1}$ such that $\mu_k(f h)=0$. In other words, we are interested in the image of $\mu_k(f)$ in $\mathcal{L}_{k+3}(H) / \mu_k\left(A_{k+1}\right)$, which is also $\mathcal{L}_{k+3}(H) / \ev_\omega\left(\tau_{k+1}'\left(A_{k+1}\right)\right)=\mathcal{L}_{k+3}(H) / \ev_\omega\left(\operatorname{Im}\left(\tau'_{k+1}\right)\right)$. 

\begin{proposition}\label{prop33}
The map $\calG_k / A_{k+1} \longrightarrow \mathcal{L}_{k+3}(H) / \ev_\omega\left(\operatorname{Im}\left(\tau'_{k+1}\right)\right)$  induced by $\mu_k$ vanishes on the image of the map $J_k/J_{k+1} \longrightarrow \calG_k/A_{k+1}$.
\end{proposition}

\begin{proof}
This follows directly from the fact that elements of $\calM$ acts trivially on $\zeta$.
\end{proof}

We now need to understand $\ev_\omega\left(\operatorname{Im}\left(\tau'_{k+1}\right)\right)$. We know that $\ev_\omega\left(\operatorname{Im}\left(\tau'_{k+1}\right) \right) \subset \ev_\omega\left( \Ker(\Tr_{k+1})\right)$, and that the equality holds stably (see Theorem \ref{satoh}). Let us denote $\overline{C_k(H)}$ the image of $D_k(H)$ by the trace $\Tr_k$. By definition we get the following commutative diagram of short exact sequences:

\begin{equation}
\label{diagpsi}
\begin{tikzcd} 
0  \arrow[r]  & \Ker(\Tr_k) \cap D_k(H) \arrow[r] \arrow[d] &D_k(H)  \arrow[r, "\Tr_k"] \arrow[d] & \overline{C_k(H)} \arrow[r] \arrow[d] & 0 \\ 0  \arrow[r] & \Ker(\Tr_k) \arrow[r] & \Der_k(H) \arrow[r, "\Tr_k"] & C_k(H)\arrow[r] & 0.
\end{tikzcd}
\end{equation}

\begin{proposition}
\label{defpsi}
There is a short exact sequence 
\[
\begin{tikzcd} 
0  \arrow[r]  & \ev_\omega\big(\Ker(\Tr_k)\big) \arrow[r] & \calL_{k+2}(H) \arrow[r, "\psi_k"] & C_k(H) /\overline{C_k(H)} \arrow[r]  & 0
\end{tikzcd} \] where for $h \in \calL_{k+1}(H)$, $\psi_k$ sends $[a_i,h]$ to the class of $\Tr_k(b_i^* \otimes h)$ and $[b_i,h]$ to the class of $\Tr_k(-a_i^* \otimes h)$.
\end{proposition}

\begin{proof}
It follows from the snake Lemma applied to the commutative diagram \eqref{diagpsi} and from the short exact sequence \eqref{ses}. To compute $\psi_k$, one can use the lifts defined in the proof of Proposition \ref{propses}.
\end{proof}

It is known that $\Ker(\Tr_k) \cap D_k(H)$ admits a structure of $\Sp(H)$-module for which $\tau_k : J_k \rightarrow \Ker(\Tr_k) \cap D_k(H)$ is equivariant over $\calM \rightarrow \Sp(H)$ (where $\calM$ acts on $J_k$ by conjugation). Also the natural action of $\Sp(H)$ on $C_k(H)$ preserves $\overline{C_k(H)}$, because $\Tr$ is $\Sp(H)$-equivariant (and $D_k(H)$ is preserved by the action of $\Sp(H)$). Hence the target space of $\psi_{k}$. For the sake of clarity, we recall that the action of $f_* \in \GL (H)$ on a derivation $d$ is $f_* \cdot d := f_* \circ d \circ f_*^{-1}$.

\begin{lemma}\label{psiequ}
The map $\psi_k : \calL_{k+2} \rightarrow C_{k}(H)/\overline{C_{k}(H)} $ is $\Sp(H)$-equivariant.
\end{lemma}

\begin{proof}
Let $f_* \in \Sp(H)$ and $l \in \calL_{k+2}$. By definition, $\psi_k(l) = \Tr_k(d)$ for any derivation $d \in \Der_k$ such that $d(\omega) = l$. Then $$(f_*\cdot d)(\omega) = f_*\cdot\big(d(f_*^{-1} \cdot \omega)\big) = f_* \cdot d(\omega)$$ hence $\psi_k(f_* \cdot l) = \Tr_k(f_* \cdot d) = f_*(\Tr_k(d))$.
\end{proof}

Summing up, we have proven the following.

\begin{theorem}
\label{thmdef}
Let $k \geq 2$. Then the map
\begin{align*}
\overline{\Tr}_k : \im(\tau'_k) \cap D_k(H) & \longrightarrow  C_{k+1}(H)/\overline{C_{k+1}(H)} \\
d & \longmapsto \psi_{k+1} \circ \mu_k(f)
\end{align*} where $f$ is any element of $A_k$ such that $\tau'_k(f) = d$ is well-defined. Moreover, it is a $\Sp(H)$-equivariant homomorphism that vanishes on the image of $\tau_k$.
\end{theorem}

\begin{proof}
The map is well-defined because by definition $f \in \calG_k$, and a different choice of lift would make the result differ by an element in $$\psi_{k+1}(\mu_k(A_{k+1})) \subset  \psi_{k+1}(\ev_\omega\left(\operatorname{Im}\left(\tau'_{k+1}\right)\right)) \subset \psi_{k+1}(\ev_\omega\left( \Ker(\Tr_{k+1})\right)) = 0.$$ It is a homomorphism by Proposition \ref{propmumorphism}, and it vanishes on $\im(\tau_k)$ by Proposition \ref{prop33}. The equivariance follows from Propositions \ref{propeq1} and \ref{psiequ}.
\end{proof}

\begin{remark}\label{deg1}
For $k = 1$, the fact that $D_1(H) = \tau'_1(J_1) = \im(\tau_1)$ implies that $\overline{\Tr}_1$ is trivial.
\end{remark}
\begin{remark}
Let $k \geq 2$ and $n = 2g \geq k+2$. Then by Theorem \ref{satoh}, the map $\overline{\Tr}_k$ can be defined on $\Ker(\Tr_k) \cap D_k(H)$.
\end{remark}

We also want to know how the map $\overline{\Tr}$ behaves with respect to brackets of derivations. The bracket of two derivations in $D(H)$ is still in $D(H)$, and the bracket of two derivations in $\im(\tau')$ is still in $\im(\tau')$, because $\tau'$ is a Lie morphism. Hence $\im(\tau') \cap D(H)$ is a subalgebra of $\Der(H)$. Any derivation acting on $H$ also acts on $\calL(H)$. We need the following proposition.

\begin{proposition}
Let $f \in \mathcal{G}_k$ and $g \in \mathcal{G}_l$ with $k, l \geq 1$. Then $[f, g] \in \mathcal{G}_{k+l}$ and:
$$
\mu_{k+l}([f, g])=\tau'_k(f)(\mu_l(g))-\tau'_l(g)(\mu_k(f)).
$$\end{proposition}

\begin{proof} First $[f,g]\in A_{k+l}$ (see Section \ref{sec2}). Then we need to compute $[f, g](\zeta) \zeta^{-1}$ modulo $\Gamma_{k+l+4}\pi$. Let us introduce some notations first:
\begin{align*}
a &:=  f(\zeta)\zeta^{-1} \in \Gamma_{k+3}\pi \\
a' &:=  f^{-1}(\zeta)\zeta^{-1} \in \Gamma_{k+3}\pi\\
b &:=  g(\zeta)\zeta^{-1} \in \Gamma_{l+3}\pi\\
b' &:=  g^{-1}(\zeta)\zeta^{-1} \in \Gamma_{l+3}\pi 
\end{align*} Remark that $a$ and $a^{\prime}$ (resp. $b$ and $b^{\prime}$ ) are linked by $f\left(a^{\prime}\right)=a^{-1}$ (resp. $g\left(b^{\prime}\right)=b^{-1}$). Now, using the relation $f(\zeta)=a \zeta$ and the other three analogous relations, we get:
$$
\begin{aligned}
{[f, g](\zeta) \zeta^{-1} } & =f g f^{-1} g^{-1}(\zeta) \zeta^{-1}=f g f^{-1}\left(b^{\prime}\right) f g\left(a^{\prime}\right) f(b) a \zeta \zeta^{-1} \\
& =f g\left(f^{-1}\left(b^{\prime}\right) b^{\prime-1}\right) f g\left(b^{\prime}\right) f\left(g\left(a^{\prime}\right) a^{\prime-1}\right) f\left(a^{\prime}\right) f(b) a \\
& =f g\left(f^{-1}\left(b^{\prime}\right) b^{\prime-1}\right) f\left(b^{-1}\right) f\left(g\left(a^{\prime}\right) a^{\prime-1}\right) f(b) f(b)^{-1} a^{-1} f(b) a \\
& =f g\left(f^{-1}\left(b^{\prime}\right) b^{\prime-1}\right) f(b)^{-1} f\left(g\left(a^{\prime}\right) a^{\prime-1}\right) f(b)\left[f(b)^{-1}, a^{-1}\right] .
\end{aligned}
$$

Then recall that $a \in \Gamma_{k+3}\pi$ and $b \in \Gamma_{l+3}\pi$, and that $f, g \in A_1$ act trivially on the quotients $\Gamma_i \pi / \Gamma_{i+1} \pi$. Thus $g\left(b^{\prime}\right)=b^{-1}$ implies that $b^{\prime}=-b$ in $\Gamma_{l+3}\pi / \Gamma_{l+4}\pi$. By applying $\tau\left(f^{-1}\right)=-\tau(f)$ (from $\mathcal{L}_{k+3}$ to $\left.\mathcal{L}_{k+l+3}\right)$ to this equality, we see that $f^{-1}\left(b^{\prime}\right) b^{\prime-1}$ belongs to $\Gamma_{k+l+3} \pi$ and is equal to $f(b) b^{-1}$ modulo $\Gamma_{k+l+4}\pi$. The same reasoning says that $g\left(a^{\prime}\right) a^{\prime-1}$ belongs to $\Gamma_{k+l+3}\pi$ and is equal to $-g(a) a^{-1}$ modulo $\Gamma_{k+l+4}\pi$. Finally, we notice that $\left[f(b)^{-1}, a^{-1}\right] \in \Gamma_{k+l+6}\pi$ is trivial modulo $\Gamma_{k+l+4}\pi$, and that conjugation by any element (in particular by $f(b)$ ) is trivial on the quotients $\Gamma_i\pi / \Gamma_{i+1}\pi$ (in particular on $\left.\Gamma_{k+l+3}\pi / \Gamma_{k+l+4}\pi\right)$. As a consequence, we have $[f, g](\zeta) \zeta^{-1} \in \Gamma_{k+l+3}\pi$ (which means that $[f, g] \in \mathcal{G}_{k+l}$ ) and we get:
$$
[f, g](\zeta) \zeta^{-1}=f(b) b^{-1}-g(a) a^{-1} \in \Gamma_{k+l+3}\pi / \Gamma_{k+l+4}\pi .
$$
\end{proof}

\begin{corollary}
The graded space $\bigoplus_{k \geq 1}\calG_k / \calG_{k+1}$ is a graded Lie ring and the map $$\mu := \oplus \mu_k : \bigoplus_{k \geq 1} \calG_k / \calG_{k+1} \rightarrow \calL_{\bullet + 3}(H)$$ is a 1-cocycle with respect of the action of $\calG$ on $\calL(H)$ given by $\tau'_k$.
\end{corollary}

Note that $\im(\tau') \cap D(H)$ acts on $C(H)$ and preserves the space $\overline{C(H)}$, because if $\Tr(d_1)= 0$ (which will be the case if $d_1 \in \im(\tau')\cap D(H)$) and $d_2 \in D(H)$, we have by Proposition \ref{cocyle}: $$ d_1 \cdot \Tr(d_2) = d_1 \cdot \Tr(d_2) - d_2 \cdot \Tr(d_1) = \Tr([d_1,d_2])\in \overline{C(H)}$$ as $[d_1,d_2] \in D(H)$. Hence  $\im(\tau') \cap D(H)$ acts on $C(H)/\overline{C(H)}$. The argument fails in degree $1$, actually, but by the same computation, using Remark \ref{remarkcocycle}, we see that $\im(\tau'_1) \cap D_1(H) = D_1(H)$ acts on $C(H)/\overline{C(H)}$.

\begin{proposition}
\label{mycocycle}
The map $\overline{\Tr}$ is a 1-cocycle with respect to the action of $\im(\tau') \cap D(H)$ on $C(H) /\overline{C(H)}$, i.e. for any $d_1, d_2 \in \im(\tau') \cap D(H)$ of degrees greater than $1$: $$ \overline{\Tr}([d_1,d_2]) = d_1 \cdot \overline{\Tr}(d_2) - d_2 \cdot \overline{\Tr}(d_1).$$
\end{proposition}

\noindent To prove the proposition, we need the following Lemma. It complements Lemma \ref{psiequ}. 

\begin{lemma}
\label{lemmapsi2}
Let $\partial \in \im(\tau'_k) \cap D_k(H)$ for $k \geq 1$, and let $l \in \calL_{\geq 3}(H)$. Then $$\psi\big(\partial(l)\big)= \partial\cdot \psi(l).$$ 
\end{lemma}

\begin{proof}
Let $d$ be any integral derivation such that $d(\omega)=l$. Then by definition, $\psi(l)= \Tr(d)$. Notice that $$[\partial,d](\omega) = \partial(d(\omega)) - d(\partial(\omega)) = \partial(l) - 0 = \partial(l).$$ Hence $\psi\big(\partial(l)\big) = \Tr([\partial,d])$ because $[\partial,d] \in \Der(H)$. Either $\Tr(\partial) = 0$ (if $k > 1$), or $\Tr(\partial) \in H$ (if $k =1$), and in both cases $d\cdot \Tr(\partial) = 0$ (as $d(H)\subset \calL(H) \subset [T(H), T(H)]$). Hence $$\psi\big(\partial(l)\big) = \Tr([\partial,d]) = \partial \cdot \Tr(d)$$ by Proposition \ref{cocyle}, which concludes.
\end{proof}

\begin{proof}[Proof of Propositon \ref{mycocycle}]
We choose two lifts $f \in A_k$ and $g \in A_l$ of $d_1$ and $d_2$, with $k \geq 2$ and $l\geq 2$. We have $\tau'_{k+l}([f,g]) = [d_1,d_2]$, so that

\begin{align*}
\overline{\Tr}([d_1, d_2]) &= \psi_{k+l+1}(\mu_{k+l}([f,g])) \\
&= \psi_{k+l+1}\big(\tau'_k(f)(\mu_l(g))-\tau'_l(g)(\mu_k(f))\big) \\
&= \tau'_k(f)\cdot \psi_{l+1}(\mu_l(g)) - \tau'_l(g)\cdot \psi_{k+1}(\mu_k(f)) \\
&= d_1 \cdot \overline{\Tr}(d_2) - d_2 \cdot \overline{\Tr}(d_1)
\end{align*} where the third equality follows from Lemma \ref{lemmapsi2}.
\end{proof}

By Remark \ref{deg1}, we obtain the following corollary.

\begin{corollary}
The map $\overline{\Tr}$ vanishes on the subalgebra of $D(H)$ generated in degree $1$.
\end{corollary}

\begin{remark}
We shall prove later, in Proposition \ref{torsprop}, that $\overline{\Tr}$ is valued in the torsion subspace of $C(H)/\overline{C(H)}$. But this can actually be deduced from the previous corollary and a result of Hain \cite{haininf} implying that $\im(\tau)^\Q$ is generated in degree $1$. Our proof does not make use of this theorem of Hain, though.
\end{remark}

\section{Computing the obstruction map}\label{sec4}

In order to compute the map $\overline{\Tr}$, we need, for a map $f$, to be able to compute $\psi(\mu(f))$, \emph{i.e.} to compute $\Tr(d)$ for any derivation $d$ such  that $d(\omega) = \mu(f)$. The infinitesimal Andreadakis-Johnson representation will provide us with such an instance of $d$.

\subsection{The infinitesimal Andreadakis-Johnson representation}\label{infajr}
For the following definitions, more details are given in \cite{masinf}. For an abelian group $M$, we denote by $M^\Q$ the rationalization $M \otimes Q$ of $M$, and if $M$ is free we view $M$ as a lattice in $M^\Q$. For example, we view $\Der(H)$ as a lattice in $\Der(H)^\Q$, which coincides with the space of derivations of $\calL^\Q$. We also denote by $\widehat{T}(H)$ (resp. $\widehat{\calL}(H)$) the degree completion of $T(H)^\Q$ (resp. $\calL (H)^\Q$). An \emph{expansion} is a monoid map $\theta : F_n \rightarrow \widehat{T}(H)$ such that $\theta(x) = 1 + \{ x \} + (\text{deg} \geq 2)$, where $\{ x \}$ mean the class of $x$ in $H$. Note in particular that $\theta$ induces a group morphism to the group $T(H)^\times$ of units of $T(H)$. The expansion is \emph{group-like} if $\theta$ is valued in the subgroup of $T(H)^\times$ of group-like elements of $\widehat{T}(H)$. Let $I$ be the augmentation ideal of $\Q[F_n]$ and $\widehat{\Q[F_n]}$ be the $I$-adic completion of the group algebra of $F_n$. A group-like expansion is actually equivalent to the choice of a complete Hopf algebra isomorphism between $\widehat{\Q[F_n]}$ and $\widehat{T}(H)$ which is the identity at the graded level. We will still denote by $\theta$ the isomorphism:

\[ \theta : \widehat{\Q[F_n]} \xrightarrow{\simeq} \widehat{T}(H). \]

An element $f \in \operatorname{Aut}(F_n)$ acts via continuous Hopf algebra automorphisms on the completion $\widehat{\Q[F_n]}$, and hence on $\widehat{T}(H)$ via $\theta$. This actions preserves the set of primitive elements of $\widehat{T}(H)$ which is naturally identified with $\widehat{\calL} (H)$. Hence, we get a map \[\rho^\theta : \operatorname{Aut}(F_n) \longrightarrow \operatorname{Aut}\big(\widehat{\calL}(H)\big).\] Furthermore, $\rho^\theta$ maps $A_1 \subset \Aut (F_n)$ to $\IAut \widehat{\calL}(H)$, \emph{i.e.} the subgroup of automorphisms that induce the identity at the graded level. Any such automorphism admits a logarithm \[ \loga (f) := \sum_{k=1}^{\infty} (-1)^{k+1}\frac{(f-Id)^{k}}{k} \] which is an element of $\Der(\widehat{\calL}(H))$ (the space of derivations of $\widehat{\calL}(H)$), which coincides with the degree completion of $\Der(H)^\Q$.

In summary, we obtain the following map: 
\begin{definition}The map
\begin{align*} 
r^\theta: A_1 & \longrightarrow \Der(\widehat{\calL}(H)) \\
f & \longmapsto \loga(\theta \circ f \circ \theta^{-1})
\end{align*}is called the \emph{infinitesimal Andreadakis-Johnson. In the formula, we still use $f$ to denote the action of $f$ on $\widehat{\Q[F_n]}$}.
\end{definition}
\noindent It has the following properties:

\begin{itemize}
\item If $f,h \in A_1$ then $r^\theta(fh) = r^\theta(f) \star r^\theta(h).$
\item If $f \in A_k$ then $r^\theta(f) = \tau_k'(f) + (deg \geq k+1)$.
\item If $f,h \in A_1$ then $r^\theta(fhf^{-1}) = e^{[r^\theta(f), -]}r^\theta(h).$
\end{itemize}

\noindent Here, for two derivations $d_1$ and $d_2$, $d_1 \star d_2$ denotes the BCH product of $d_1$ and $d_2$: 

\[ d_1 \star d_2 := d_1 + d_2 + \frac{1}{2}[d_1,d_2] + \dots\]

\begin{notation}
Since the space $\Der(\widehat{\calL}(H))$ is graded, for any derivation $d$, we denote $d_k$ the part of degree $k$. Note that $\Der_i(\widehat{\calL}(H)) \simeq H^* \otimes \widehat{\calL}_{i+1}(H).$ Then the second property above can be reformulated as: for $f$ in $A_k$, $r_i^\theta(f) = 0$ if $i < k$ and $r_k^\theta(f) = \tau'_k(f).$ 
\end{notation}

\subsection{The infinitesimal Dehn-Nielsen representation}

The representation $r^\theta$ defined in Section \ref{infajr} also restricts to the Torelli group. In the case of a surface, we like to use \emph{symplectic expansions}. A symplectic expansion is a group-like expansion (see Section \ref{infajr}) satisfying the additional property that $\theta(\zeta) = e^{-\omega}$. It is proven in \cite[Lemma 2.16]{masinf} that such expansions do exist.

From now on $\theta$ is a symplectic expansion, and the \emph{infinitesimal Dehn-Nielsen representation} is defined to be the restriction of $r^\theta$ to the Torelli group $\calI = J_1$.

\begin{proposition}
\label{propr_symp}
For any $f \in \mathcal{I}$, we have that $r^\theta(f) \in D(H)^\Q \subset \Der(H)^\Q$.
\end{proposition}

\begin{proof}
It suffices to compute $r^\theta(f)(\omega) = \loga(\rho^\theta(f))(\omega)$. As $\theta(\loga(\zeta^{-1})) = \loga( e^{\omega}) = \omega$ we notice that 

\begin{align*}
(\rho^\theta(f) - Id)(\omega) &= \theta \Big(f\big(\loga(\zeta^{-1})\big)\Big) -\omega \\
&= \theta\Big(\loga\big(f(\zeta^{-1})\big)\Big) - \omega \\
&= \theta\big(\loga(\zeta^{-1})\big) - \omega. \\
&= 0
\end{align*}
It is easy to conclude  that $r^\theta(f)(\omega) = 0$.
\end{proof}

In the next section, we will also need the following lemma, which is of course coherent with Proposition \ref{propr_symp}. Let us denote $l^\theta := \loga\circ~ \theta.$

\begin{proposition}
\label{importantlemma}
Let $k \geq 2$, let $f$ be in $A_k$, and let $\theta$ be a symplectic expansion. Then \[r^\theta(f)(\omega) = \tau_k(f)(\omega) -l_{k+3}^\theta\big(f(\zeta)\zeta^{-1}\big) + (deg \geq k+4).\] Equivalently $r_{k+1}^\theta(f)(\omega) = -\mu_k(f)$.
\end{proposition}

\begin{proof}
As $(\rho^\theta(f) - Id)$ shifts the degree by $k \geq 2$, we have that, up to elements of degree greater than $k+4$, $r^\theta(f)(\omega) = (\rho^\theta(f) - Id)(\omega)$. As in the proof of Proposition \ref{propr_symp}, we get 

\begin{align*}
\big(\rho^\theta(f) - Id\big)(\omega) &=  \theta\Big(\loga\big(f(\zeta^{-1})\big)\Big) - \omega  \\
&= - l^{\theta}\big(f(\zeta)\big) - \omega \\
&= - l^{\theta}\big(f(\zeta)\zeta^{-1}) \big) \star l^\theta(\zeta) - \omega.
\end{align*}
As $f(\zeta) \zeta^{-1} \in \Gamma_{k+2}\pi$, its image by $l^\theta$ starts in degree $k+2$, whereas $l^\theta(\zeta)$ starts in degree $2$. Hence, by the BCH formula, up to degree $k+4$ we get \begin{align*}
\big(\rho^\theta(f) - Id\big)(\omega) &\equiv - \Big(l^{\theta}\big(f(\zeta)\zeta^{-1} \big) + l^\theta(\zeta)\Big) - \omega \\ 
&\equiv -l^{\theta}\big(f(\zeta)\zeta^{-1})
\end{align*} because $\theta$ is symplectic. We conclude that $r_k^\theta(f)(\omega) = -l_{k+2}^\theta(f(\zeta)\zeta^{-1})= \tau_k(f)(\omega)$ and that $r_{k+1}^\theta(f)(\omega) = -l_{k+3}^\theta(f(\zeta)\zeta^{-1})$.
\end{proof}

This gives us the following fundamental fact.

\begin{proposition}
For any $d \in \im(\tau'_k) \cap D_k(H)$, we have $\overline{\Tr}_k(d) = -\psi_{k+1}(r_{k+1}^\theta(f)(\omega))$.
\end{proposition}

\noindent The following proposition generalizes, \cite[Prop. 5.3]{massak}, quoted above as Proposition \ref{massak}:

\begin{proposition}
\label{massakimproved}
For any $k \geq 2$, for any $n \geq 2$, for any $f \in A_k$ and for any group-like expansion $\theta$ we have
$\Tr^\Q(r^\theta_{k+1}(f))= 0.$
\end{proposition}

\begin{proof}
The proof follows from the proof of \cite[Prop. 5.3]{massak}. It is shown, using a theorem of Higman and Stallings about the Whitehead group of a free group and the so-called Magnus representation, that the reduction of $\Tr^\Q(r^\theta(f))$ modulo the two-sided ideal generated by $\calL^\Q_{\geq k+1}(H)$ must be 0. In the reference, this corresponds to the statement ``Therefore $\Tr^{\mathfrak{R}}\operatorname{Ab}^\mathfrak{R}(\psi)$ starts in degree $k\geq 2$, so that it must be $0$", where $\mathfrak{R}$ is the two-sided ideal in question and $\psi$ is $r^\theta(f)$. Notice that, in degree $k+1$, the ideal $\mathfrak{R}$ is trivial since an element of $\calL^\Q_{k+1}$ is trivial in $C_{k+1}^\Q$ (it will always be a linear combination of elements of the form $[x,h] = xh-hx = 0 \in C^\Q(H)$). Hence we get the result. 
\end{proof}

\section{Properties of the infinitesimal Dehn-Nilsen representation}\label{sec5}

In this section, we investigate some properties of $r^\theta$, and deduce some facts about the map $\overline{\Tr}$, defined in Theorem \ref{thmdef}. Recall that for any group-like expansion $\theta$, we set $l^\theta := \loga \circ~\theta$. We now enter in a series of technical results that are of major importance for the proofs of Theorem A and B.

\begin{proposition}
\label{propcentral}
For any $k \geq 2$, for any $f \in A_k$, for any $x \in \pi$ we have  
\[r^\theta_{k+1}(f)(\{x\}) = l_{k+2}^\theta\big(f(x)x^{-1}\big) + \frac{1}{2}[\tau'_k(f)(\{ x\}),\{x\}] - \tau'_k(f)\big(l_2^\theta(x)\big) \in \calL_{k+2}^\Q(H). \]
\end{proposition}

\begin{proof}
We use the fact, as $k \geq 2$, that $$l^\theta_{k+1,k+2}(f(x)x^{-1}) = \theta_{k+1,k+2}(f(x)x^{-1})$$ because $f(x)x^{-1} \in \Gamma_{k+1} \pi$. The latter also implies that $l^\theta_{k+1}(f(x)x^{-1})$ is the class of $f(x)x^{-1}$ in $\calL_{k+1}(H)$, hence is equal, by definition, to $\tau'_k(f)(\{ x\})$. We then consider the following computation, where $\equiv$ stands for equivalence modulo degree greater than $k+3$:

\begin{align*}
\theta\big( f(x) \big) - \theta(x)&=  \theta\big( f(x)x^{-1} \big) \theta(x) - \theta(x) \notag \\ 
&\equiv \Big(1+ \theta_{k+1,k+2}\big( f(x)x^{-1} \big) \Big)\theta(x) - \theta(x) \notag \\
&\equiv l^\theta_{k+1}(f(x)x^{-1})(1+\{x\}) + l_{k+2}^\theta(f(x)x^{-1}) \notag\\
&\equiv \tau'_k(f)(\{ x\})(1+\{x\}) + l_{k+2}^\theta(f(x)x^{-1}).
\end{align*}

\noindent The map $\rho^\theta(f) - Id$ shifts the degree by $k$, since $f \in A_k$. Then $r^\theta(f)(\{ x\}) :=\loga(\rho^\theta(f))(\{ x\})$ is equivalent to $\big(\rho^{\theta}(f)-Id\big) (\{ x\})$ modulo terms of degree $2k+1 \geq k+3$.  

\begin{align*}
r^\theta_{k+1,k+2}(f)(\{ x\})&\equiv  \big(\rho^{\theta}(f)-Id\big)(\{ x\}) \\ 
&\equiv \Big(\rho^{\theta}(f)-Id\Big)\big(\{ x\} - \theta(x)\big) +  \rho^{\theta}(f)\big(\theta(x)\big) - \theta(x)\\
&\equiv \Big(\rho^{\theta}(f)-Id\Big)\big( - \theta_{\geq 2}(x)\big) + \theta\big( f(x) \big) - \theta(x).
\end{align*}
\noindent Then, $\theta_{\geq 3}$ is shifted to degree greater than $k+3$, and we have that $\theta_2(x) = l_2^\theta(x) + \frac{\{ x\}^2}{2}$. Hence $$\Big(\rho^{\theta}(f)-Id\Big)\big(  \theta_{2}(x)\big) \equiv \tau'_k(f)\big(l_2^\theta(x) + \frac{\{ x\}^2}{2}\big).$$ Hence, by the computation of $\theta\big(f(x)\big)$ made above, we get: 

\begin{align*}
r^\theta_{k+1,k+2}(f)(\{ x\})&\equiv  \Big(\rho^{\theta}(f)-Id\Big)\big( - \theta_{2}(x)\big) + \theta\big( f(x) \big) - \theta(x) \\
&\equiv -\tau'_k(f)\big(l_2^\theta(x) + \frac{\{ x\}^2}{2}\big) + \tau'_k(f)(\{ x\})(1+\{x\}) + l_{k+2}^\theta(f(x)x^{-1}) \\
&\equiv \tau'_k(f)(\{ x\}) + l_{k+2}^\theta(f(x)x^{-1}) +\frac{1}{2}[\tau'_k(f)(\{ x\}),\{x\}] - \tau'_k(f)\big(l_2^\theta(x)\big).
\end{align*} where $\tau'_k(f)(\{ x\}^2)$ means applying the unique extension of $\tau'_k(f)$ to a derivation of $T(H)$ to $ \{ x\}^2$, so that $$\tau'_k(f)(\{ x\}^2)= \tau'_k(f)(\{ x\})\{ x\}  + \{ x\}\tau'_k(f)(\{ x\}).$$
\end{proof}

For the next proposition, we need to define an operation $\delta$ that associates to a parenthesized word with letters in the basis of $H$ an element of $\calL(H)$.
\begin{definition} Consider a parenthesized word with $k \geq 2$ letters in the basis of $H$ (we shall omit the external parentheses). To such an element, we can associate the following  element of $\calL_{k+1}(H)$. Build the planar rooted binary tree associated to the parenthesized word and orient it from root to leaves. Then take the sum over all the leaves of the Lie algebra element associated to the tree obtained by doubling the leaf with the following sign rule: 

\[ 
\raisebox{-7mm}{\includegraphics[scale=1.35]{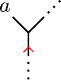}}\hspace{4mm} \mapsto \hspace{4mm} \raisebox{-7mm}{\includegraphics[scale=1.35]{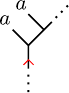}} \text{     and     } 
\raisebox{-7mm}{\includegraphics[scale=1.35]{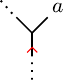}}\hspace{4mm} \mapsto \hspace{4mm}- \raisebox{-7mm}{\includegraphics[scale=1.35]{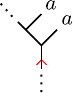}}
\]
\end{definition}
\begin{remark}
As an operation from colored  rooted binary trees to $\calL(H)$, the operation $\delta$ passes to the quotient by the antisymmetry relation and the Jacobi relation. However, it does not behave well if we allow the letters to be any elements of $H$.
\end{remark}

\begin{example}
If $u,v,w$ are letters in the basis of $H$, we have $$\delta( uv) = [u,[u,v]] - [[u,v],v]$$ and $$\delta\big( u(vw)\big) = [u,[u,[v,w]]] + [u,[v,[v,w]]] - [u,[[v,w],w]]. $$
\end{example}

\noindent We set the following notation. 
\begin{notation}Denote by $d_{l^\theta_2}$ the unique derivation of $\calL^\Q(H)$ which sends $\{ x \}$ to $l_2^\theta(x)$ for $x$ in the chosen basis for $\pi$.
\end{notation}

\noindent The operation $\delta$ is interesting to us because of the following lemma.

\begin{lemma}
\label{deltalemma}
Let $k \geq 2$, and $\tilde{\gamma} \in \Gamma_k \pi$ be a commutator of elements of the chosen basis for $\pi$, represented by a parenthesized word $w_\gamma$ of length $k$. Let $\gamma$ be the class of $\tilde{\gamma}$ in $\calL_{k}(H)$. Then $$l_{k+1}^\theta(\tilde{\gamma}) = d_{l_2^\theta}(\gamma) + \frac{1}{2}\delta(w_\gamma).$$
\end{lemma}

\begin{proof}
We prove the lemma by induction on $k \geq 2$. We simply need the following formula  for the BCH commutator:
\begin{align}
\label{eqbch}
[x,y]_{\star} := x \star y \star (-x) \star (-y) = [x,y] + \frac{1}{2} [x,[x,y]] - \frac{1}{2}[y,[y,x]] + \text{commutators of length $\geq 4$}. 
\end{align} If $k = 2$, then $l_3^\theta([u,v])$ is obtained as the degree 3 part of the BCH commutator of $l^\theta(u)$ and $l^\theta(v)$. Also for $u \in \Gamma_k \pi$, $l^\theta_k(u)$ is simply the class of $u$ in $\Gamma_k \pi / \Gamma_{k+1} \pi \simeq \calL_{k}(H)$. hence we get 
\begin{align*}
l^\theta_3([u,v]) &= [l_1^\theta(u),l_2^\theta(v)] + [l_2^\theta(u),l_1^\theta(v)] + \frac{1}{2}[l_1^\theta(u),[l_1^\theta(u),l_1^\theta(v)]] - \frac{1}{2}[l_1^\theta(v),[l_1^\theta(v),l_1^\theta(u)]] \\
&= [\{u\},l_2^\theta(v)] + [l_2^\theta(u),\{v\}] + \frac{1}{2}[\{u\},[\{u\},\{v\}]] - \frac{1}{2}[\{v\},[\{v \},\{u\}]]  \\
&= d_{l_2^\theta}([\{u\}, \{v \}]) + \frac{1}{2} \delta(uv).
\end{align*}
Now, if $\tilde{\gamma}$ has commutator length $k+1$, then it is a commutator of elements of smaller sizes. If it is a commutator of two elements $\tilde{\gamma_1}$ and $\tilde{\gamma_2}$ of length $r_1, r_2$ greater than $2$, then the BCH formula for the commutator yields, in degree $k+2$: 
\[ l^\theta_{k+2}([\tilde{\gamma_1}, \tilde{\gamma_2}]) = [l^\theta_{r_1}(\tilde{\gamma_1}), l^\theta_{r_2+1}(\tilde{\gamma_2})] + [l^\theta_{r_1+1}(\tilde{\gamma_1}),l^\theta_{r_2}(\tilde{\gamma_2})]\] and we can conclude by induction. Indeed, it is clear from the definition that $\delta((w_1)(w_2))= [\delta(w_1),w_2] + [w_1, \delta(w_2)]$, for two words $w_1$ and $w_2$ of length greater than $2$ (note that there is abuse of notation here, as we considered $w_1$ and $w_2$ as elements in $\calL(H)$). If, however, $\tilde{\gamma} = [\tilde{\gamma_1}, x]$ for $x$ in the basis of $\pi$, then we conclude by using the BCH formula as in the case $k = 2$.

\end{proof}
We extend $\delta$ linearly to linear combination of parenthesized words in the basis of $H$. Let $f$ be in $A_k$. For all $x$ in $H$, we choose a decomposition of $\tau_k'(f)(x)$ as a sum of commutators of elements of the basis of $H$. Wee see this is as a sum of parenthesized words, and we denote (\emph{abusively, as it depends on the choice of decomposition}) by $\delta(\tau_k'(f)(x))$ its image by $\delta$. We set, for $x$ ranging in the basis of $H$, $$\delta(\tau'_k(f)) := \sum_{x} x^* \otimes \delta(\tau'_k(f)(x)).$$
\begin{lemma}
\label{technical}
Let $k \geq 2$, and $f \in A_k$, then the class of $r^\theta_{k+1}(f)$ in $\Der(H)^\Q $ modulo the lattice of integral derivations is given by the following formula:
\[ r^\theta_{k+1}(f) \equiv  [d_{l_2^\theta}, \tau_k'(f)] + \frac{1}{2}\delta(\tau_k'(f)) + \frac{1}{2}\sum_{x } x^* \otimes [\tau_k'(f)(x), x] \text{ } \operatorname{mod} \text{ } \Der_{k+1}(H)  \] where the sum ranges over the basis of $H$, and the space of derivations is once again identified with $\operatorname{Hom}(H, \calL(H)) = H^* \otimes \calL(H)$.
\end{lemma}

\begin{proof}
For each $x$ in the basis of $H$, consider the product $\gamma_x$ of the ``obvious" lifts of the elements in the chosen decomposition of $\tau'_k(f)(x)$. Let $\tilde{x}$ be the lift of $x$ in the basis of $\pi$. For example $[x,[y,z]] - [y,[y,x]]$ would be lifted to $[\tilde{x}, [\tilde{y}, \tilde{z}]] [\tilde{y},[\tilde{y},\tilde{x}]]^{-1}$. Then by definition $e_x := f(\tilde{x})\tilde{x}^{-1}\gamma_{x}^{-1} \in \Gamma_{k+2}\pi$. We now apply Proposition \ref{propcentral} to get the following:

\begin{align*}
r^\theta_{k+1}(f)(x) &= l_{k+2}^\theta\big(e_x\gamma_x\big) + \frac{1}{2}[\tau_k(f)( x),\{x\}] - \tau_k(f)\big(l_2^\theta(\tilde{x})\big) \\
&= l^\theta_{k+2}(e_x) + l_{k+2}^\theta(\gamma_x) + \frac{1}{2}[\tau_k(f)( x),\{x\}] - \tau_k(f)\big(l_2^\theta(\tilde{x})\big),
\end{align*}
since $l^\theta_{k+2}$ is additive on $\Gamma_{k+1}\pi$. Because of the same additivity property, Lemma \ref{deltalemma} then implies that $$l_{k+2}^\theta(\gamma_x) = d_{l_2^\theta}\big(\tau_k'(f)(x)\big) +\frac{1}{2} \delta(\tau'_k(f)(x)).$$ Also $e_x$ is in $\Gamma_{k+2}\pi$, hence its image by $l_{k+2}^\theta$ is integral.
We deduce, in $\calL_{k+2}(H)^\Q/\calL_{k+2}(H)$:

\begin{align*}
r^\theta_{k+1}(f)(x) &\equiv l_{k+2}^\theta\big(\gamma_x\big) + \frac{1}{2}[\tau'_k(f)( x),x] - \tau'_k(f)\big(l_2^\theta(\tilde{x})\big) \\
&\equiv  d_{l_2^\theta}(\tau_k'(f)(x)) + \frac{1}{2}\delta(\tau'_k(f)(x)) + \frac{1}{2}[\tau'_k(f)( x),x] - \tau'_k(f)\big(d_{l_2^\theta}(x)\big).
\end{align*} The conclusion follows.
\end{proof}

We can now deduce the following proposition. 

\begin{proposition}
\label{torsprop}
The map $\overline{\Tr}$ is valued in the $2$-torsion subspace of $C(H)/\overline{C(H)}$.
\end{proposition}

\begin{proof}
We have shown that $\overline{\Tr}(d)$ can be computed as $-\psi_{k+1}(r^\theta_{k+1}(f)(\omega))$ for $\tau'_k(f) = d$. Hence we can choose the symplectic expansion from \cite{masinf}, so that $l_2^\theta$ is valued in $\frac{1}{2}\calL_{2}(H) \subset \calL^\Q_{2}(H)$. This implies that $d_{l^\theta_2} \in \frac{1}{2}\Der_{1}(H)$, hence we deduce from Lemma \ref{technical} that for $f \in A_k$, $r^\theta_{k+1}(f) \in \frac{1}{2}\Der_{k+1}(H)$. Furthermore, we know from Proposition \ref{massakimproved} that $\Tr^\Q(r^\theta_{k+1}(f)) = 0$. Hence we clearly have \[ r^\theta_{k+1}(f^2) = 2r^\theta_{k+1}(f) \in \Ker(\Tr).\] Consider an element $d \in \im(\tau'_k) \cap D_k(H)$, and a lift $f$ of $d$ to $A_k$. Then $f^2$ is a lift of $2d$, hence $\overline{\Tr}_k(2d) = \psi_{k+1}\big( r^\theta_{k+1}(f^2)(\omega) \big)= \Tr_{k+1}(r^\theta_{k+1}(f^2))= 0$.
\end{proof}

\section{Tools for computing the trace}\label{sec6}

Before computing the map $\overline{\Tr}$ on some examples, we would like to have a better understanding of its target space $C(H)/\overline{C(H)}$.

\subsection{Target space of the obstruction}

It is noted by Conant \cite{conant} that the map $\Tr$, when restricted to $D(H)$, has an  ``extra $\mathbb{Z}_2$ symmetry". We clarify this point by studying $\overline{C(H)}$. First, we recall a useful diagrammatic description due to Levine. More details can be found in \cite{faesmas}. A \emph{tree diagram} is a finite, unitrivalent, connected, acyclic graph whose trivalent vertices are oriented (i.e$.$ edges are cyclically ordered around each trivalent vertex), and whose univalent vertices are colored by~$H$. The orientation in the tree diagrams in this paper is given by the trigonometric orientation (all the diagrams will be drawn in the plane). The degree of a tree diagram is the number of trivalent vertices in the tree. A \emph{leaf} of the tree is a univalent vertex. The \emph{space of tree diagrams} $\calT(H)$ is the abelian group generated by tree diagrams modulo the relations in Figure \ref{fig1}.

\begin{figure}[h]
\vspace{12pt}
\centering
\begin{subfigure}[b]{0.32\textwidth}
\centering
${\tikz[baseline=10pt, line width = 0.25mm, scale = 0.6]{
\draw[color=wqwqwq] (0.75,0)-- (0.75,0.6);
\draw[color=wqwqwq] (0.75,0.6) arc (-80:70:0.3) -- (0,1.5);
\draw[color=wqwqwq] (0.75,0.6) arc (-100:-240:0.3);
\draw[color=wqwqwq] (0.77,1.25)-- (1.5,1.5)}} = - {\tikz[baseline=10pt, scale = 0.6, line width = 0.25mm,]{
\draw[color=wqwqwq] (0.75,0)-- (0.75,0.75);
\draw[color=wqwqwq] (0.75,0.75) -- (0,1.5);
\draw[color=wqwqwq] (0.75,0.75) --(1.5,1.5) }}$

 \caption*{$\text{AS}$}
\end{subfigure}
\hfill
\begin{subfigure}[b]{0.32\textwidth}
\centering
${\tikz[baseline=10pt, line width = 0.25mm,  scale = 0.6]{
\draw [color=wqwqwq] (0,0)-- (1.5,0);
\draw [color=wqwqwq] (0.75,0)-- (0.75,1.5);
\draw [color=wqwqwq] (0,1.5)-- (1.5,1.5);}} = {\tikz[baseline=10pt,line width = 0.25mm, scale = 0.6]{
\draw [color=wqwqwq] (0,0)-- (0,1.5);
\draw [color=wqwqwq] (0,0.75)-- (1.5,0.75);
\draw [color=wqwqwq] (1.5,0)-- (1.5,1.5);}} - {\tikz[baseline=10pt,line width = 0.25mm, scale = 0.6]{
\draw [color=wqwqwq] (0,0)-- (0.7,0.7);
\draw [color=wqwqwq] (0.8,0.8)-- (1.5,1.5);
\draw [color=wqwqwq] (0,1.5)-- (1.5,0);
\draw [color=wqwqwq] (0.35,0.35)-- (1.15,0.35);}}$
\caption*{IHX}
\end{subfigure}
\hfill
\begin{subfigure}[b]{0.32\textwidth}
\centering
$
{\tikz[baseline=10pt, scale = 0.6, line width = 0.25mm,]{
\draw[color=wqwqwq] (0.75,0)-- (0.75,0.75);
\draw[color=wqwqwq] (0.75,0.75) -- (0,1.5);
\draw[color=wqwqwq] (0.75,0.75) --(1.5,1.5);
\node[label= {[xshift=0cm, yshift=0.75cm]{$a+b$}}]{} }} = {\tikz[baseline=10pt, scale = 0.6, line width = 0.25mm,]{
\draw[color=wqwqwq] (0.75,0)-- (0.75,0.75);
\draw[color=wqwqwq] (0.75,0.75) -- (0,1.5);
\draw[color=wqwqwq] (0.75,0.75) --(1.5,1.5);
\node[label= {[xshift=0cm, yshift=0.75cm]{$a$}}]{} }} + 
{\tikz[baseline=10pt, scale = 0.6, line width = 0.25mm,]{
\draw[color=wqwqwq] (0.75,0)-- (0.75,0.75);
\draw[color=wqwqwq] (0.75,0.75) -- (0,1.5);
\draw[color=wqwqwq] (0.75,0.75) --(1.5,1.5);
\node[label= {[xshift=0cm, yshift=0.75cm]{$b$}}]{} }}$

 \caption*{multilinearity}
\end{subfigure}

\caption{The $\text{AS}$, $\text{IHX}$, and multilinearity relations}
\label{fig1}
\end{figure}

The space $\calT(H)$ is a graded Lie ring where the bracket of two trees $S$ and $T$ is obtained by summing all the ways of $\omega$-connecting a leaf of $S$ with a leaf of $T$. By \emph{$\omega$-connecting} two vertices, we mean that we glue together the vertices, forget the labels, and multiply by the contraction of the labels by $\omega : H \otimes H \rightarrow \Z$ (see \cite{faesmas}).

For a leaf $v$ in a tree $T$ we denote $T_v$ the rooted tree obtained by forgetting the label $v$. We denote by $\operatorname{lie}(T_v)$ the element of $\calL(H)$ induced by the tree $T_v$. For example 
$$\operatorname{lie}\Bigg(\begin{tikzpicture}[baseline =-10mm, x=0.75pt,y=0.75pt,yscale=-0.6,xscale=0.6]

\draw    (20,20) -- (80,80) ;
\draw    (140,20) -- (80,80) ;
\draw    (80,20) -- (110,50) ;
\draw    (110,20) -- (95,35) ;
\draw    (80,80) -- (80,104) ;

\draw (17,8) node {$a$};
\draw (77,5) node {$b$};
\draw (107,8) node  {$c$};
\draw (137,5) node   {$d$};

\end{tikzpicture}
\Bigg) =  [a,[[b,c],d]] \in \calL_4(H).
\quad
$$ \vspace{0.1cm}

\noindent The \emph{expansion} map $\eta : \calT(H) \rightarrow D(H)$ is obtained by sending a tree $T$ to \[\eta(T) := \sum_{v} \omega(v,-)\otimes \operatorname{lie}(T_v)\] where the sum is taken over all the leaves $v$ of $T$.

It is known that, rationally, the map $\eta$ is an isomorphism. Using the work of Levine \cite{levtrees} and Conant, Schneiderman and Teichner \cite{CST}, we claim the following more precise facts.

\begin{proposition}[Levine, Conant-Schneiderman-Teichner]
\label{LCST}
Let $k \geq 1$ be an integer. 
\begin{itemize}
\item $D_{2k}(H)$ is isomorphic, via the inverse of $\eta^\Q$, to the sublattice of $\calT_{2k}(H)^\Q$ generated by $\calT_{2k}(H)$ and the elements $\frac{1}{2} u \trait u$, where $u$ is any rooted tree of degree $k$.
\item $D_{2k+1}(H)$ is isomorphic, via the inverse of $\eta^\Q$, to the quotient of $\calT_{2k+1}(H)$ by the relation in Figure \ref{fig2}, for all $x \in H$ and $T$ rooted tree of degree $k+1$.
\end{itemize} 
\end{proposition}

\begin{figure}[h]

$$
{\tikz[baseline=10pt, line width = 0.25mm,]{
\draw[color=wqwqwq] (0.75,0)-- (0.75,0.75);
\draw[color=wqwqwq] (0.75,0.75) -- (0,1.5);
\draw[color=wqwqwq] (0.75,0.75) --(1.5,1.5);
\node[label= {[xshift=0cm, yshift=1.35cm]{$T$}}]{} ;
\node[label= {[xshift=1.5cm, yshift=1.35cm]{$T$}}]{} ;
\node[label= {[xshift=0.75cm, yshift=-0.5cm]{$x$}}]{} ;} } = 0. $$
\caption{A relation for $D_{2k+1}(H)$}
\label{fig2}
\end{figure}

\begin{remark}
By the $AS$ relation, the tree on the left-side in Figure \ref{fig2} is $2$-torsion in $\calT_{2k+1}(H)$. Actually, $D_{2k+1}(H)$ can be identified with the image of $\calT_{2k+1}(H)$ in $\calT_{2k+1}(H)^\Q$. 
\end{remark}

Now, following \cite[Section 4]{conant}, we recall the diagrammatic interpretation of the Satoh trace. Recall that any element in $H^* \otimes \calL(H) = \operatorname{Hom}(H, \calL(H))$ corresponds bijectively to the restriction of a derivation of $\calL(H)$. Similarly to the description made above, if we draw a rooted tree $T$ with an element $h^* $ of $H^*$ at the root, we mean the derivation $h^* \otimes \operatorname{lie}(T)$. Let us call a \emph{necklace} a diagram consisting of one oriented cycle on which are glued rooted $H$-colored tree diagrams. The space of necklaces $\calN(H)$ is the abelian group generated by necklaces subject to the IHX and AS relations (even at the trivalent vertices on the cycle) and to the multilinearity relation. It is graded by the number of trivalent vertices. As explained in the proof of \cite[Theorem 4.2]{conant}, one can read a necklace by running along the cycle, and by reading the cyclic word obtained by using the IHX relations whenever one meets a trivalent vertex (see the example below). Notice that one should read positively a letter on the left-side. Then we have 

\begin{itemize}
\item If $h^* \otimes T \in \Der_k(H)$ then $\Tr(h^* \otimes T)$ is the element of $C_k(H)$ obtained by summing over all the way of contracting $h^*$ with a leaf of $T$ and reading the necklace obtained by adding an oriented edge towards the leaf.
\item If $T \in D(H)$ then $\Tr(T)$ coincides with the element of $C_k(H)$ obtained by summing over all the ways of $\omega$-connecting two leaves $a$ and $b$ of $T$, and reading the necklace obtained by adding an oriented edge from $a$ to $b$. 
\end{itemize}

\begin{example}
We clarify the discussion above by showing an example. Let $a,b,c,d,e$ be different basis elements in $H$. Then

\tikzset{every picture/.style={line width=0.75pt}} 
$$ \Tr \Bigg(
\begin{tikzpicture}[baseline = -15mm, x=0.75pt,y=0.75pt,yscale=-0.7,xscale=0.7]

\draw    (121,39.11) -- (60.1,112.78) ;
\draw    (16,82.56) -- (60.1,112.78) ;
\draw    (68.5,56.11) -- (89.5,76.89) ;
\draw    (26.5,54.22) -- (68.5,56.11) ;
\draw    (70.6,24) -- (68.5,56.11) ;
\draw    (60.1,112.78) -- (60,136.6) ;
\draw    (45,32.8) -- (47.5,55.17) ;

\draw (2,68.4) node [anchor=north west][inner sep=0.75pt]    {$a$};
\draw (16,42.4) node [anchor=north west][inner sep=0.75pt]    {$b$};
\draw (40,12.4) node [anchor=north west][inner sep=0.75pt]    {$c$};
\draw (65,4.4) node [anchor=north west][inner sep=0.75pt]    {$d$};
\draw (55,134.4) node [anchor=north west][inner sep=0.75pt]    {$d^{*}$};
\draw (118,19.4) node [anchor=north west][inner sep=0.75pt]    {$e$};
\end{tikzpicture} \Bigg) = 
\begin{tikzpicture}[baseline = -15mm, x=0.75pt,y=0.75pt,yscale=-0.7,xscale=0.7]

\draw    (136,54.31) -- (75.1,127.98) ;
\draw    (31,97.76) -- (75.1,127.98) ;
\draw    (84.55,70.37) -- (105.55,91.14) ;
\draw    (41.5,69.42) -- (83.5,71.31) ;
\draw    (85.6,39.2) -- (83.5,71.31) ;
\draw    (75.1,127.98) -- (75,151.8) ;
\draw  [color={rgb, 255:red, 208; green, 2; blue, 27 }  ,draw opacity=1 ] (152.91,97.41) -- (159.75,88.77) -- (162.56,99.42) ;
\draw    (60,48) -- (62.5,70.37) ;
\draw [color={rgb, 255:red, 208; green, 2; blue, 27 }  ,draw opacity=1 ]   (85.6,39.2) .. controls (88,-0.4) and (159,11.6) .. (161,44.6) .. controls (163,77.6) and (160,99.6) .. (148,127.6) .. controls (136,155.6) and (73,182.6) .. (75,151.8) ;

\draw (17,83.6) node [anchor=north west][inner sep=0.75pt]    {$a$};
\draw (31,57.6) node [anchor=north west][inner sep=0.75pt]    {$b$};
\draw (132,33.6) node [anchor=north west][inner sep=0.75pt]    {$e$};
\draw (55,29.6) node [anchor=north west][inner sep=0.75pt]    {$c$};

\end{tikzpicture}
 =        
 \begin{tikzpicture}[baseline = -15mm, x=0.75pt,y=0.75pt,yscale=-0.7,xscale=0.7]

\draw  [color={rgb, 255:red, 208; green, 2; blue, 27 }  ,draw opacity=1 ] (11.6,81.02) .. controls (11.6,60.31) and (28.39,43.52) .. (49.1,43.52) .. controls (69.81,43.52) and (86.6,60.31) .. (86.6,81.02) .. controls (86.6,101.73) and (69.81,118.52) .. (49.1,118.52) .. controls (28.39,118.52) and (11.6,101.73) .. (11.6,81.02) -- cycle ;
\draw  [color={rgb, 255:red, 208; green, 2; blue, 27 }  ,draw opacity=1 ] (45.38,113.02) -- (55.27,117.87) -- (45.45,122.87) ;
\draw    (49.1,43.52) -- (48.6,64.52) ;
\draw    (72.6,32.52) -- (65.6,47.72) ;
\draw    (89.6,41.52) -- (77.6,56.72) ;
\draw    (25.6,33.52) -- (28.6,49.72) ;
\draw  [color={rgb, 255:red, 208; green, 2; blue, 27 }  ,draw opacity=1 ] (132.6,81.02) .. controls (132.6,60.31) and (149.39,43.52) .. (170.1,43.52) .. controls (190.81,43.52) and (207.6,60.31) .. (207.6,81.02) .. controls (207.6,101.73) and (190.81,118.52) .. (170.1,118.52) .. controls (149.39,118.52) and (132.6,101.73) .. (132.6,81.02) -- cycle ;
\draw  [color={rgb, 255:red, 208; green, 2; blue, 27 }  ,draw opacity=1 ] (166.38,113.02) -- (176.27,117.87) -- (166.45,122.87) ;
\draw    (170.1,43.52) -- (169.6,64.52) ;
\draw    (193.6,32.52) -- (186.6,47.72) ;
\draw    (210.6,41.52) -- (198.6,56.72) ;
\draw    (146.6,33.52) -- (149.6,49.72) ;
\draw    (102,78.4) -- (116,78.4) ;

\draw (44.6,62.92) node [anchor=north west][inner sep=0.75pt]    {$e$};
\draw (19.6,16.92) node [anchor=north west][inner sep=0.75pt]    {$a$};
\draw (68.6,15.92) node [anchor=north west][inner sep=0.75pt]    {$b$};
\draw (88.6,21.92) node [anchor=north west][inner sep=0.75pt]    {$c$};
\draw (165.6,62.92) node [anchor=north west][inner sep=0.75pt]    {$e$};
\draw (140.6,16.92) node [anchor=north west][inner sep=0.75pt]    {$a$};
\draw (209.6,21.92) node [anchor=north west][inner sep=0.75pt]    {$b$};
\draw (189.6,13.92) node [anchor=north west][inner sep=0.75pt]    {$c$};

\end{tikzpicture} = -cbea + bcea.
$$
Indeed, as one can check: $\Tr(d^* \otimes [a,[[[b,c],d],e]]) = d^*(d)([b,c]ea) = -bcea +cbea.$
\end{example}

\subsection{Mirrors}
When computing the trace of an element in $D(H)$, one can see from the graphical interpretation, that for each unordered pair of leaves $\{a,b \}$, we get a term $\omega(a,b)w$ and a term $\omega(b,a)(-1)^{k}\overline{w}$, where $w$ is a sum of words, and $\overline{w}$ the mirror image of $w$.

\begin{definition} For $k \geq 1$, we denote $\operatorname{Mir}_k(H)$ the subspace of $C_k(H)$ generated by elements of the form $w + (-1)^{k+1}\overline{w}$. \end{definition}
\noindent Recall that, by definition, $\overline{C_k(H)} = \Tr_k(D_k(H))$. As we just explained, \cite[Theorem 4.2]{conant} and its proof imply
\begin{proposition}
For any $k \geq 1$, $\overline{C_k(H)} \subset \operatorname{Mir}_k(H)$.
\end{proposition}

\noindent The equality holds in degree $2$, since $\operatorname{Mir}_2(H) = 0$, and the trace vanishes on $D_2(H)$. The equality also seems to hold in degree 3, hence the following question:

\begin{question}
\label{question}
Do we have, at least for $n \gg k$, $\overline{C_k(H)} = \operatorname{Mir}_k(H) $ ?
\end{question}

\noindent According to \cite[Corollary 5.6]{conant}, rationally, the answer to Question \ref{question} is yes, i.e. the quotient $ \operatorname{Mir}_k / \overline{C_k(H)}$ is torsion.

We now define a reduction of $\overline{\Tr}$:

\begin{definition}
\begin{equation}
\label{deftrmir}
\begin{tikzcd}
\Tr_k^{\mir} : \im(\tau_k') \cap D_k(H) \arrow[r, "\overline{\Tr_k}"] & C_{k+1}(H)/\overline{C_{k+1}(H)}  \arrow[r] & C_{k+1}(H)/\operatorname{Mir}_{k+1}(H).
\end{tikzcd}
\end{equation}
\end{definition} 

\noindent We switch to the (simpler) study of $C_k(H)/\operatorname{Mir}_k(H)$. It was pointed out to the author, after sharing a first version of this paper, that the next result was already proven in \cite[Prop. 5.2]{nss}. We still give a proof of this result, as we shall later need some elements in the proof. In the sequel, a \emph{chiral} necklace, or chiral tensor, is a cyclic tensor equal to its mirror image (up to cyclic permutation).

\begin{proposition}
\label{bkspace}
Let $\varphi$ be Euler's totient function. We have the following description of $B_k(H) := C_k(H)/\operatorname{Mir}_k(H)$:

\begin{itemize}
\item If $k$ is even, then $B_k(H)$ is free abelian of rank $$\frac{1}{2k}\displaystyle\sum_{d \mid k}\varphi(d)n^{k/d} + \frac{1}{4}(n+1)n^{k/2}.$$
\item If $k$ is odd, then $B_k(H)$ is the direct sum of a free abelian part of rank $$\frac{1}{2k}\displaystyle\sum_{d \mid k}\varphi(d)n^{k/d} - \frac{1}{2}n^{(k+1)/2}$$ and a 2-torsion part of rank $n^{k+1/2}.$ 
\end{itemize}
\end{proposition}

\begin{proof}
First, we notice that if we pick a basis of $H$, and consider the associated basis of $T_k(H)$, the action of $\Z/k\Z$ by cyclic permutation preserves this basis, which explains why $C_k(H)$ is free abelian of rank the number of different necklaces of size $k$ with $n$ beads. Now, we further quotient by identifying $w$ with $(-1)^{k}\overline{w}$. This action also preserves the basis, up to sign. Hence, we claim that if $k$ is even (i.e. there is no sign), $B_k(H)$ is free abelian with rank the number of \emph{bracelet} of size $k$ with $n$ beads. A bracelet is a necklace up to reflection. When $k$ is odd, a non-chiral necklace will be identified with the opposite of its mirror (dividing by $2$ the dimension of the subspace of non-chiral necklaces), while a chiral necklace will be identified to minus itself (introducing 2-torsion with rank the number of chiral necklaces). Hence we reduce the problem to computing the number of chiral necklaces, and the number of necklaces. We denote $N_{k,n}$ and $B_{k,n}$ the number of necklaces and bracelets of size $k$ with $n$ colors of beads. We can use Burnside's lemma to compute this number, looking at the actions of $\Z / k\Z$ and the dihedral group $D_{2k}$ on words of size $k$. The difference between the odd and even cases comes from the fact that the number of conjugacy classes in $D_{2k}$ is $2$ or $1$ depending on the parity of $k$. We get the following formulas: 

\begin{align*}
N_{k,n} &= \frac{1}{k}\displaystyle\sum_{d \mid k}\varphi(d)n^{k/d}, \\
B_{k,n} &= \begin{cases}\frac{1}{2k}\displaystyle\sum_{d \mid k}\varphi(d)n^{k/d} + \frac{1}{4}(n+1)n^{k/2} & \text { if } k \text { is even } \\ \frac{1}{2k}\displaystyle\sum_{d \mid k}\varphi(d)n^{k/d} + \frac{1}{2}n^{(k+1)/2} & \text { if } k \text { is odd. }\end{cases}
\end{align*}
\noindent But $B_{k,n}$ is also the number of chiral necklaces plus half the number of non-chiral necklaces, while $N_{k,n}$ is equal to the number of chiral necklaces plus the number of non-chiral necklaces. Hence we see that the number of chiral necklaces is $\frac{1}{2}(n+1)n^{k/2}$ when $k$ is even and $n^{(k+1)/2}$ otherwise. 
\end{proof}

\begin{corollary}
The map $\Tr^{\mir}$ is trivial in odd degrees, and valued in the subspace generated by chiral necklaces in even degrees.
\end{corollary}

\begin{proof}
We have proven in Proposition \ref{bkspace} that the space $B_k(H)$ has no torsion if $k$ is even. $\Tr_{k}^{\mir}$ is valued in the 2-torsion subspace of $ B_{k+1}(H)$, hence is trivial if $k$ is odd. When $k$ is even, we know from the proof of Proposition \ref{bkspace} that the 2-torsion of $B_{k+1}$ is exactly the subspace generated by the projections of chiral necklaces.  
\end{proof}

\begin{remark}
The reduction from $C_k(H)/\overline{C_k(H)}$ to $B_k(H)$ a priori kills some torsion, hence we can not affirm that $\overline{\Tr}_k$ vanishes in odd degrees.
\end{remark}

\subsection{Some examples, and the degree 2 case}
We now give some examples, using the description of $D(H)$ given by Proposition \ref{LCST}. Note that $D_2(H) \subset \im(\tau_2')$ so that $\overline{\Tr}_2$ and $\Tr_{2}^{\mir}$  are well-defined on $D_2(H)$.
\begin{example}
\label{stupid}
The tree $\frac{1}{2}\vltree{a_1}{b_1}{a_1}{b_1}$, representing the  derivation $$b_1^* \otimes [[b_1,a_1],b_1] - a_1^*\otimes [a_1,[b_1,a_1]]\in D_2(H) \cap \Ker(\Tr),$$ is in the image of $\tau_2$. Hence we have $\overline{\Tr}(\frac{1}{2}\vltree{a_1}{b_1}{a_1}{b_1}) = 0.$ 
\end{example}

For a non trivial example, we consider the following elements of $A_1$, for $x,y, z$ some generators of $\pi$. The automorphism $K_{xy}$ sends $x$ to $y^{-1}xy$ and fixes the other generators, while $K_{xyz}$ sends $x$ to $x[y,z]$ and fixes the other generators. To lighten the notation in the computation below, for $x$ in the basis of $H$, we also use the notation $x$ for its lift in the basis of $\pi$. In the last section of this paper, we shall give a formula for computing $\overline{\Tr}$, which is why we allow ourselves not to give every detail in the computation below.

\begin{example}
\label{annoying}
The derivation $d := b_1^* \otimes [a_2,[a_1,a_2]] + b_2^* \otimes [a_1,[a_2,a_1]] \in \Ker(\Tr) \cap D_2(H)$, given in Example \ref{exampled}, admits a lift in $A_2$ given by $f := [K_{a_1a_2}, K_{b_1a_1a_2}][K_{a_2a_1}, K_{b_2a_2a_1}]$. Using a Sagemaths computer program, we get 
\begin{align*}
f(\zeta)\zeta^{-1} =& 
a_2^2a_1a_2^{-1}a_1^{-1}a_2a_1^{-1}a_2^{-1}a_1b_2^{-1}a_2^{-1}b_2a_1^{-1}a_2a_1 \\
&a_2^{-1}(a_1a_2a_1^{-1}a_2^{-1}a_1)^2a_2^{-1}a_1^{-1}a_2b_1^{-1}a_1^{-1}b_1a_2^{-1}\\&a_1a_2a_1^{-1}a_2a_1a_2^{-1}a_1^{-1}b_1^{-1}a_1b_1a_1^{-1}b_2^{-1}a_2b_2a_2^{-1}.
\end{align*}
By direct computation, we can then check that \begin{align*}
-r_3^\theta(f)(\omega) &= l^\theta_{5}(f(\zeta)\zeta^{-1}) \\
&= -[a_1,[a_1,[[a_1,a_2],a_2]]]-[a_1,[[[a_1,a_2],a_2],a_2]]+[a_1,[a_1,[a_2,[a_2,b_2]]]]\\
& \hspace{4mm}+[a_1,[[a_1,[a_2,b_2]],a_2]]-[[a_1,a_2],[[a_1,a_2],a_2]]-[[a_1,a_2],[a_1,[a_2,b_2]]]\\&\hspace{4mm}+[[[a_1,a_2],a_2],[a_1,b_1]].
\end{align*}
Then we can check directly, using Proposition \ref{defpsi}, that 
\begin{align*}
\overline{\Tr}(d) &= \psi_{3}\big([[a_1,a_2],[a_1,[a_2,b_2]]]-[[[a_1,a_2],a_2],[a_1,b_1]]\big) \\
&= \Tr(-b_2^*\otimes [a_1,[a_1,[a_2,b_2]]] - a_1^* \otimes [[[a_1,a_2],a_2],a_1]) \\
&= a_2a_1a_1 - a_2a_2a_1.
\end{align*} Hence we see that $\overline{\Tr}(d)$ is non-zero.
\end{example}

\noindent In fact, we can prove the following relation between ${\Tr}^{\mir}_2$ and $\Tr^{as}$. Let $\iota : \Lambda^2(H) \otimes \Z_2\rightarrow B_3(H)$ be defined by $\iota(x \wedge y) = xxy + yyx$. Notice that $xxy = -xxy \in B_3(H)$, because $xxy$ is a chiral tensor (equal to its mirror image up to cyclic permutation).

\begin{proposition}
\label{tracesrel}
We have the following commutative diagram of $\Sp(H)$-linear maps: 
\[
\begin{tikzcd}
D_2(H) \arrow[r, "{\Tr}^{\mir}_2"]\arrow[rd, "\Tr^{as}",bend right =20] &B_3(H)\\
& \Lambda^2(H) \otimes \Z_2. \arrow[u,"\iota"]
\end{tikzcd} \]
\end{proposition}

\begin{proof}
As all of the maps are equivariant, and both traces vanish on $\im(\tau_2)$, it is enough to prove the commutativity for lifts of ($\Sp(H)$-module) generators of the cokernel of $\tau_2$. According to \cite[Theorem 2.4]{faes}, this cokernel is isomorphic (via $\Tr^{as}$) to $ \operatorname{Ker}(\omega : \Lambda^2(H/2H) \rightarrow \mathbb{Z}_2 )$ which is generated by $a_1 \wedge a_2$. Hence it is enough to check the equality on the element $d$ in Example \ref{annoying}. The conclusion then follows from Example \ref{exampled}, since, as explained, $\Tr^{as}(d) = a_1 \wedge a_2$. 
\end{proof}

\begin{remark}
Since $\iota$ is injective, we have $\Ker(\Tr^{as}) = \Ker(\Tr_2^{\mir})$, hence $\Ker(\Tr_2^{\mir}) = \im(\tau_2)$. We deduce that the kernel of $\overline{\Tr}_2$ is also equal to $\im(\tau_2)$, and that $\overline{\Tr}$ is not richer that $\Tr^{\mir}$ in degree 2.
\end{remark}

\section{Torsion in the cokernels of the Johnson homomorphisms}
\label{sec7}

As one can see from Example \ref{annoying}, it is not so easy to compute directly the map $\overline{\Tr}$. In particular, finding a lift in $A_k$ of an element in  $\Ker(\Tr_k)$ can be quite challenging. Even the rest of the computation can hardly be done by hand. We shall now provide a formula that allows us to compute the map on a derivation given under the form of a sum of \emph{elementary trees}.

\subsection{A closed formula for the computation of $\overline{\Tr}$}

\begin{definition} By an elementary tree diagram we mean a tree diagram colored by elements of the basis.By a ``half of a symmetric tree", we mean an element of shape  $\frac{1}{2} u \trait u$, where $u$ is any rooted tree of degree $k$. \end{definition}
 
Recall that an element of $D(H)$ will be represented by a sum of elements of $\calT(H)$ and (in the even cases) halves of symmetric trees (see Proposition \ref{LCST}).  
 
\begin{definition}[Doubling trees]
For an elementary tree $T$, we set $\Delta(T) \in \calT_{k+1}(H) $ to be the \emph{double} of $T$, i.e. the sum of elementary trees obtained by doubling each external edge once in the following way: if $T$ locally looks like $\hspace{1mm}\raisebox{-4mm}{\includegraphics[scale=1]{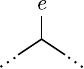}}\hspace{1mm}$ then the tree obtained by doubling the edge colored by $e$ should look like $\hspace{1mm}\raisebox{-3mm}{\includegraphics[scale=1]{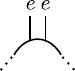}}\hspace{1mm}$. For a ``half of symmetric tree" $\frac{1}{2} u \trait u$, we set $\Delta(T) =\frac{1}{2}\Delta(u \trait u)$, and this element is still in $\calT_{k+1}(H)$ because of the symmetries.

\end{definition}
The operation is similar to $\delta$ except that there is no way to chose an orientation in this case. This operation does not respect the antisymmetry or the multilinearity, and is \emph{not} an operation on the space $\calT(H)$. 

\begin{example}
For $a,b,c \in H$, we have $$\Delta(\ltritree{a}{b}{c}) = \vltree{a}{b}{c}{a} +  \vltree{b}{c}{c}{a} +  \vltree{a}{b}{b}{c}$$ and $$\Delta(\frac{1}{2}\vltree{a}{b}{a}{b}) = \vlfivetree{a}{b}{a}{b}{a} +\vlfivetree{a}{b}{a}{b}{b}.$$ 
\end{example}

\begin{remark}
The operation $\Delta$ can be extended to the vector space generated by tree diagrams, but it does not respect the antisymmetry or the multilinearity relations, hence it does  \emph{not} define an operation on the space $\calT(H)$. 
\end{remark}

\begin{notation}For an elementary tree diagram $T$, we set $E(T)$ the set of external vertices of $T$. For $t \in E(T)$ a colored external vertex of $T$, we denote by $T_t$ the rooted tree obtained by forgetting $t$ and turning in into the root. We also set $\overline{a_i}:= b_i$, $\overline{b_i} := a_i$, and $x^\omega := \omega(x,-)$. Notice that $a_i^\omega = \omega(a_i, -) = \overline{a_i}^*$ and $b_i^\omega =\omega(b_i, -) = -\overline{b_i}^*$.\end{notation} 

\noindent We shall now use tree notations, as they are much easier to read: this is meant to help the reader during the course of the proof of the next proposition, which is somewhat challenging.

\begin{proposition}
\label{prophard}
Let $d$ be an element of $\im(\tau') \cap D(H)$, written as a $\Z$-linear combination of elementary trees and possibly halves of symmetric elementary  trees. Let $\calS$ be the set of these trees, and $\lambda_T$ the integral coefficient in front of $T$. We have \[ \overline{\Tr}(d) = \sum_{T \in \calS} \lambda_T \Bigg( \sum_{t \in E(T)}\Big(\Tr\big(\rtritree{t}{T_t}{t^*}\big)+ \frac{1}{2} \Tr\big(\rtritree{t + \overline{t}}{T_t}{\overline{t}^*}\big)\Big) - \frac{1}{2}\Tr\big(\Delta(T)\big) \Bigg).\]
\end{proposition}

\begin{remark}
One could think that this formula extends $\overline{\Tr}$ to all the elements of $D(H)$. This is not the case, as there is no reason, in general, for the result to be independent of the choice of the decomposition in elementary trees.
\end{remark}

\begin{proof}[Proof of Proposition \ref{prophard}]
We have proven in Lemma \ref{technical}, that for any choice of $f$ in $A_k$ we have \[ r^\theta_{k+1}(f) = X + [d_{l_2^\theta}, \tau_k'(f)] + \frac{1}{2}\delta(\tau_k'(f)) + \frac{1}{2} \sum_{x } x^* \otimes [\tau_k'(f)(x), x] \] for some element $X \in \Der_{k+1}(H)$. Here, if $f$ is a lift of $d$ (\emph{i.e.} if $\tau'_k(f) = d$), $\delta(\tau'_k(f))$ is given by \[ \delta(\tau_k'(f)):= \sum_{T \in \calS} \lambda_T \sum_{t \in E(T)} \overline{t}^* \otimes \delta(T_t). \] where $T_t$ is seen as a parenthesized words in $H$.
We use the following observation: since $X$ is an integral derivation, $\psi_{k+1}(X(\omega)) = \Tr(X)$ by definition. Also, the trace of $r^\theta_{k+1}(f)$ is trivial, so we have $$\Tr\Big([d_{l_2^\theta}, \tau_k'(f)] + \frac{1}{2}\delta(\tau_k'(f)) + \frac{1}{2}\sum_{x } x^* \otimes [\tau_k'(f)(x), x])\Big) = \Tr(-X) \in C_k(H).$$ We set $$F:= \frac{1}{2}\delta(\tau_k'(f)) + \frac{1}{2}\sum_{x } x^* \otimes [\tau_k'(f)(x), x]) \in \frac{1}{2} H^* \otimes L_{k+2}(H) = \frac{1}{2}\Der_{k+1}(H).$$ We remark that $\Tr([d_{l_2^\theta}, \tau_k'(f)]) = 0$, because of Remark \ref{remarkcocycle}, hence $\Tr(F) = - \Tr(X)$. Moreover, since $\tau'_k(f) = d \in D_k(H)$, we have $\tau'_k(f)(\omega) = 0$, so we get  \begin{align*}
\overline{\Tr}(d) &= \psi_{k+1}\big(r^\theta_{k+1}(f)(\omega)\big) \\
&= \psi_{k+1}\big(X(\omega)\big) + \psi_{k+1}\big( [d_{l_2^\theta}, \tau_k'(f)](\omega) + F(\omega) \big) \\
&= \Tr(X) + \psi_{k+1}\big( [d_{l_2^\theta}, \tau_k'(f)](\omega) + F(\omega) \big)  \\
&= \psi_{k+1}\Big( -\tau_k'(f)\big(d_{l^\theta_2}(\omega)\big) + F(\omega) \Big) - \Tr(F). 
\end{align*}
We want to emphasize that the latter term can be computed without lifting the derivation $d$ to an element of $A_k$, since it only depends on $\tau'_k(f) = d$. We used the definition of $\psi$ to get rid of the term in $X$. We shall now use a second trick to get rid of the most complicated term, the one involving the operation $\delta$. We set $$G := \frac{1}{2}\sum_{T \in \calS} \lambda_T \sum_{t \in E(T)} \rtritree{T_t}{\overline{t} - t}{t^\omega}.$$ We prove in Lemma \ref{lemmadouble} below that there exists an integral derivation $Y$ such that $$F - \frac{1}{2} \sum_{T\in \calS} \lambda_T\Delta(T) = Y + G. $$ It is clear that $\Delta$ takes value in $\calT(H)$ which is identified with the space of symplectic derivations, so for every $T \in \calS$, $\Delta(T)(\omega) = 0.$ Also, $\psi_{k+1}(Y(\omega)) =  \Tr(Y)$, because $Y$ is integral. From this we deduce \begin{align*}
\overline{\Tr}(d) &= \psi_{k+1}\Big( -\tau_k'(f)\big(d_{l^\theta_2}(\omega)\big) + F(\omega) \Big) - \Tr(F) \\ 
&= \psi_{k+1}\Big( -\tau_k'(f)\big(d_{l^\theta_2}(\omega)\big) + G(\omega) + Y(\omega)\Big) - \Tr(F) \\
&= \psi_{k+1}\Big( -\tau_k'(f)\big(d_{l^\theta_2}(\omega)\big) + G(\omega)\Big) + \Tr(Y) - \Tr\big(Y+G +\frac{1}{2} \sum_{T\in \calS} \lambda_T \Delta(T)\big) \\
&= \psi_{k+1}\Big( -\tau_k'(f)\big(d_{l^\theta_2}(\omega)\big) + G(\omega)\Big) - \Tr\big(G +\frac{1}{2} \sum_{T\in \calS} \lambda_T \Delta(T)\big).
\end{align*}  
We shall now compute this term, by direct computation, using the special value of $\theta$ given in \cite{masinf}. For this particular $\theta$, we have \begin{align*}
d_{l_2^\theta}(a_i) &= - \frac{1}{2}[a_i,b_i] \\
d_{l_2^\theta}(b_i) &= - \frac{1}{2}[a_i,b_i].
\end{align*}
We shall use rooted trees to write elements in the free Lie algebra. We have $$-d_{l_2^\theta}(\omega)= \frac{1}{2}\sum_{j = 1}^g \treethreel{a_j}{b_j}{b_j}  + \treethreer{a_j}{a_j}{b_j}.
$$ We first compute \begin{align*}
-\tau_k'(f)\big(d_{l^\theta_2}(\omega)\big) = \frac{1}{2}\sum_{T \in \calS} \lambda_T \sum_{t \in E(T)} \sum_{j = 1}^g &t^\omega(a_j)\bigg(\treethreer{T_t}{a_j}{b_j}+ \treethreer{a_j}{T_t}{b_j} +\treethreel{T_t}{b_j}{b_j}\bigg)\\
 +& t^\omega(b_j)\bigg(\treethreer{a_j}{a_j}{T_t} + \treethreel{a_j}{T_t}{b_j}+  \treethreel{a_j}{b_j}{T_t}\bigg)
\end{align*} and 
\begin{align*}
G(\omega)& = \frac{1}{2}\sum_{T \in \calS} \lambda_T \sum_{t \in E(T)} \sum_{j = 1}^g t^\omega(a_j)\btreethreel{T_t}{\overline{t} - t}{b_j} + t^\omega(b_j)\btreethreer{a_j}{T_t}{\overline{t} - t}. \\
\end{align*} We then notice that $t^\omega(a_j) \neq 0$ if and only if $t=b_j$ and $t^\omega(b_j) \neq 0$ if and only if $t= a_j$ (where $t =x$ is an abuse of language to say that $t$ is labeled by $x$). Hence, when computing the sum of the two terms above, for $t = b_j$, we get terms of the form \begin{align*}
&t^\omega(a_j)\bigg(\treethreer{T_t}{a_j}{b_j}+ \treethreer{a_j}{T_t}{b_j} +\treethreel{T_t}{b_j}{b_j} + \treethreel{T_t}{a_j}{b_j} - \treethreel{T_t}{b_j}{b_j}\bigg) \\
=& t^\omega(a_j)\bigg(2\treethreel{T_t}{a_j}{b_j} -  2\treethreel{T_t}{b_j}{a_j}\bigg)
\end{align*} and for $t = a_j$, terms of the form

\begin{align*}
&t^\omega(b_j)\bigg(\treethreer{a_j}{a_j}{T_t} + \treethreel{a_j}{T_t}{b_j}+  \treethreel{a_j}{b_j}{T_t} + \treethreer{a_j}{T_t}{b_j} - \treethreer{a_j}{T_t}{a_j} \bigg) \\
=& t^\omega(b_j)\bigg(2\treethreel{a_j}{T_t}{b_j} + 2\treethreer{a_j}{a_j}{T_t} \bigg).
\end{align*} This means that we have

$$-\tau_k'(f)\big(d_{l^\theta_2}(\omega)\big) +G(\omega) = \sum_{T \in \calS} \lambda_T \sum_{t \in E(T)} \sum_{j = 1}^g t^\omega(a_j)\bigg(\treethreel{T_t}{a_j}{b_j} -  \treethreel{T_t}{b_j}{a_j}\bigg) + t^\omega(b_j)\bigg(\treethreel{a_j}{T_t}{b_j} + \treethreer{a_j}{a_j}{T_t} \bigg), $$ which is, as it should be, an element of $\calL_{k+3}(H)\subset \calL^\Q_{k+3}(H).$ Indeed, if $T$ is integral, then the term in the sum is integral, while if $T$ is half a symmetric tree, then each term appear twice (once for each copy of $t$). We can now compute the image by $\psi_{k+1}$, using Proposition \ref{defpsi}:
\begin{align*}
\psi_{k+1}\Big( -\tau_k'(f)\big(d_{l^\theta_2}(\omega)\big) + G(\omega)\Big)=
\sum_{T \in \calS} \lambda_T \sum_{t \in E(T)} \sum_{j = 1}^g &t^\omega(a_j)\Tr\bigg( \rtritree{T_t}{a_j}{a_j^*}+\rtritree{T_t}{b_j}{b_j^*}\bigg)\\ +& t^\omega(b_j)\Tr\bigg( \rtritree{a_j}{T_t}{a_j^*}+\rtritree{a_j}{T_t}{b_j^*}\bigg)
\end{align*} Note that this formula is correct even when considering halves of symmetric trees, as one can check by grouping the terms $T_t$ for symmetric copies of $t$. We then use the fact that $t^\omega = t^\omega(\overline{t})\overline{t}^*= \pm \overline{t}^* $, where $t^\omega(\overline{t})$ is equal to~$1$ if $t = a_i$ and $-1$ if $t= b_i$, and compute:
\begin{align*}
\Tr(G) &= \frac{1}{2} \sum_{T \in \calS} \lambda_T \Tr\bigg( \rtritree{T_t}{\overline{t}}{t^\omega} - \rtritree{T_t}{t}{t^\omega}\bigg) \\ 
&= \frac{1}{2} \sum_{T \in \calS} \lambda_T \sum_{t \in E(T)} \sum_{j = 1}^g t^\omega(a_j)\Tr\bigg( \rtritree{T_t}{a_j}{a_j^*} - \rtritree{T_t}{b_j}{a_j^*} \bigg) + t^\omega(b_j)\Tr\bigg( \rtritree{T_t}{b_j}{b_j^*} - \rtritree{T_t}{a_j}{b_j^*} \bigg)
\end{align*}
and $\psi_{k+1}\Big( -\tau_k'(f)\big(d_{l^\theta_2}(\omega)\big) + G(\omega)\Big) -\Tr(G)$ is equal to

\begin{align*}
\sum_{T \in \calS} \lambda_T \sum_{t \in E(T)} \sum_{j = 1}^g &t^\omega(a_j)\Tr\bigg( \frac{1}{2}\rtritree{T_t}{a_j}{a_j^*}+\rtritree{T_t}{b_j}{b_j*} + \frac{1}{2} \rtritree{T_t}{b_j}{a_j*}\bigg) \\&+ t^\omega(b_j)\Tr\bigg( \rtritree{a_j}{T_t}{a_j*} + \frac{1}{2} \rtritree{a_j}{T_t}{b_j^*} - \frac{1}{2} \rtritree{T_t}{b_j}{b_j^*}\bigg).
\end{align*}

\noindent Finally, using $t^\omega = t^\omega(\overline{t})\overline{t}^*$ and examining separately the cases where $t= a_j$ and $t = b_j$, we get:

\begin{align*}
\overline{\Tr}(d) &= \psi_{k+1}\Big( -\tau_k'(f)\big(d_{l^\theta_2}(\omega)\big) + G(\omega)\Big) - \Tr\big(G +\frac{1}{2} \sum_{T\in \calS} \lambda_T \Delta(T)\big) \\
&= \sum_{T \in \calS} \lambda_T \Bigg(\sum_{t \in E(T)} \bigg(\Tr\Big( \rtritree{t}{T_t}{t^*}\Big) + \frac{1}{2}\Tr\Big( \rtritree{\overline{t}}{T_t}{\overline{t}^*} + \rtritree{t}{T_t}{\overline{t}^*}\Big)\bigg) - \frac{1}{2}\Tr\big(\operatorname{D(T)}\big)\Bigg)\\
&= \sum_{T \in \calS} \lambda_T \Bigg( \sum_{t \in E(T)}\bigg(\Tr\Big(\rtritree{t}{T_t}{t^*}\Big)+ \frac{1}{2} \Tr\Big(\rtritree{t + \overline{t}}{T_t}{\overline{t}^*}\Big)\bigg) - \frac{1}{2}\Tr\big(\Delta(T)\big) \Bigg).
\end{align*}

\end{proof}

\begin{lemma}
With the notation in the proof of Proposition \ref{prophard}, we have
\label{lemmadouble}
$$F - \frac{1}{2} \sum_{T\in \calS}\lambda_T \Delta(T) \equiv \frac{1}{2}\sum_{T \in \calS}\lambda_T \sum_{t \in E(T)} \rtritree{T_t}{\overline{t}-t}{t^\omega}  \operatorname{mod} \Der_{k+1}(H).$$
\end{lemma}

\begin{proof}
Let us write $L_t$ and $R_t$ for the left and right part of the tree $T_t$ (i.e $[L_t, R_t] = T_t$):

$$ T = \rtritree{L_t}{R_t}{t}.$$ Recall that the operation $\Delta$ is defined similarly to $\delta$, but without any sign convention. Choose a vertex $t$ of $T$. If we consider only the terms in the sum $\Delta(T)$ for which $t$ has not been doubled, then all the terms obtained by expanding through the expansion map $\eta$ where $t$ is taken as the root will coincide up to signs with $\omega(t,-)\delta(T_t)$. On the contrary if $t$ is the vertex that has been doubled, two terms appear:  $$\vltree{L_t}{t^\omega}{t}{R_t} +  \vltree{L_t}{t}{t^\omega}{R_t}.$$ Recall that $\tau'k(f) =d = \sum_{T \in \calS}\lambda_T T$. We have: 

\begin{align*}
\frac{1}{2}\sum_{T \in \calS}\lambda_T\Delta(T) &=  \frac{1}{2}\sum_{T \in \calS}\lambda_T\sum_{t \in E(T)} \pm t^\omega \otimes \delta(T_t) + \vltree{L_t}{t^\omega}{t}{R_t} +  \vltree{L_t}{t}{t^\omega}{R_t} \\ 
&\equiv \frac{1}{2} \delta(\tau_k'(f)) + \frac{1}{2}\sum_{T \in \calS}\lambda_T\sum_{t \in E(T)} \vltree{L_t}{t^\omega}{t}{R_t}- \vltree{R_t}{t^\omega}{t}{L_t}\\
&\equiv \frac{1}{2} \delta(\tau_k'(f)) + \frac{1}{2}\sum_{T \in \calS}\lambda_T\sum_{t \in E(T)} \rtritree{T_t}{t}{t^\omega}.
\end{align*} Here we changed the sign in front of $\vltree{R_t}{t^\omega}{t}{L_t}$, because we work in $\frac{1}{2}\Der_{k+1}(H)/ \Der_{k+1}(H)$, and then used the IHX relation. But once again, we have to be extra careful to the case where $T_t$ is not an integral tree. In this case the term appears twice, and we can change the sign in front of both terms, so that the computation is still valid. Furthermore, by direct diagrammatic interpretation, we have 
$$ \sum_{x } x^* \otimes [\tau_k'(f)(x), x]) =\sum_{T \in \calS}\lambda_T\sum_{t \in E(T)} \rtritree{T_t}{\overline{t}}{t^\omega}. $$ Hence 

$$ F-\frac{1}{2}\sum_{T \in \calS}\lambda_T\Delta(T) \equiv  \frac{1}{2}\sum_{T \in \calS}\lambda_T \sum_{t \in E(T)} \rtritree{T_t}{\overline{t}-t}{t^\omega}  \operatorname{mod} \Der_{k+1}(H).$$ 
\end{proof}

\subsection{Examples, and the proof of Theorem B}

We shall now use Proposition \ref{prophard} to compute $\overline{\Tr}$ on some examples. We divide the computation in 3 terms, $\overline{\Tr}^i$ for $i = 1,2,3$ corresponding to the 3 kinds of terms in the proposition.

\begin{align*}
\overline{\Tr}^1(d) &:= \sum_{T \in \calS} \lambda_T \sum_{t \in E(T)}\Tr\Big(\rtritree{t}{T_t}{t^*}\Big) \\
\overline{\Tr}^2(d) &:= \sum_{T \in \calS} \lambda_T \sum_{t \in E(T)} \frac{1}{2} \Tr\Big(\rtritree{t + \overline{t}}{T_t}{\overline{t}^*}\Big)\\
\overline{\Tr}^3(d) &:= \sum_{T \in \calS} \frac{\lambda_T}{2}\Tr\big(\Delta(T)\big).
\end{align*} This is purely notational, as these terms depend on some choices made (and not only on the derivation $d$), while $\overline{\Tr}(d) = \overline{\Tr}^1(d) +\overline{\Tr}^2(d) - \overline{\Tr}^3(d)$ depends only on the derivation $d$.

\begin{example}
We first check coherence with Example \ref{stupid}, we get \begin{align*}
\overline{\Tr}_2^1\big(\frac{1}{2}\vltree{a_1}{b_1}{a_1}{b_1}\big) &= a_1a_1b_1 + b_1b_1a_1 \\
\overline{\Tr}_2^2\big(\frac{1}{2}\vltree{a_1}{b_1}{a_1}{b_1}\big) &= -\frac{1}{2}\big( b_1a_1(a_1+b_1)+a_1b_1(b_1+a_1) \big) \\
&= -a_1a_1b_1 - b_1b_1a_1 \\ 
\overline{\Tr}_2^3\big(\frac{1}{2}\vltree{a_1}{b_1}{a_1}{b_1}\big) &= \Tr\big(\frac{1}{2}\Delta(\frac{1}{2}\vltree{a_1}{b_1}{a_1}{b_1})\big) \\
&= 0.
\end{align*} Hence $\overline{\Tr} \big(\frac{1}{2}\vltree{a_1}{b_1}{a_1}{b_1} \big) = 0,$ as expected.
\end{example}

\begin{example}
\label{ex513}
We check coherence with Example \ref{annoying}, we get \begin{align*}
\overline{\Tr}_2\big(\frac{1}{2}\vltree{a_1}{a_2}{a_1}{a_2}\big) &= \overline{\Tr}_2^1\big(\frac{1}{2}\vltree{a_1}{a_2}{a_1}{a_2}\big)\\
&= a_2a_2a_1 + a_1a_1a_2,
\end{align*} as expected. Indeed the word $2a_2a_1a_1$ is in $\overline{C_3(H)}$, hence $a_2a_1a_1 \equiv - a_2a_1a_1$.
\end{example}

\begin{example}
We compute $\overline{\Tr}\big( \vltree{a_1}{a_3}{a_2}{b_1} \big)$:
\begin{align*}
\overline{\Tr}_2^1\big(\vltree{a_1}{a_3}{a_2}{b_1}\big) &= 0 \\
\overline{\Tr}_2^2\big(\vltree{a_1}{a_3}{a_2}{b_1}\big) &=- \frac{1}{2} \big(a_2a_3(a_1+b_1) + a_3a_2(a_1+b_1)\big)\\ 
\overline{\Tr}_2^3\big(\vltree{a_1}{a_3}{a_2}{b_1}\big) &= a_3a_2a_2 + a_2a_3a_3 - \frac{1}{2}\big(a_3(a_1+b_1)a_2 + a_2(a_1+b_1)a_3 \big)\\
&= 0. 
\end{align*} Hence, $$\overline{\Tr}\big( \vltree{a_1}{a_3}{a_2}{b_1} \big) = -a_3a_2a_2 - a_2a_3a_3 = \iota\Big(\Tr^{as}\big( \vltree{a_1}{a_3}{a_2}{b_1} \big)\Big)$$ as expected. 
\end{example}

We now prove Theorem B, by generalizing Example \ref{ex513}:

\begin{proof}[Proof of Theorem B]
We consider the following element of $D_{2k}$, which is clearly in $\Ker(\Tr_{2k})$, and hence is stably in $\im(\tau'_{2k})$: \[ d_{(a_1,a_2, \cdots, a_{k+1})} := \frac{1}{2}\begin{tikzpicture}[baseline = -14mm, x=0.75pt,y=0.75pt,yscale=-0.7,xscale=0.7]

\draw    (50,40) -- (50,100) ;
\draw    (50,70) -- (280,70) ;
\draw    (280,40) -- (280,100) ;
\draw    (80,40) -- (80,70) ;
\draw    (110,40) -- (110,70) ;
\draw    (150,40) -- (150,70) ;
\draw    (180,100) -- (180,70) ;
\draw    (220,100) -- (220,70) ;
\draw    (250,100) -- (250,70) ;

\draw (117,47) node [anchor=north west][inner sep=0.75pt]    {$\cdots $};
\draw (187,80.4) node [anchor=north west][inner sep=0.75pt]    {$\cdots $};
\draw (42,102) node [anchor=north west][inner sep=0.75pt]    {$a_{1}$};
\draw (271,23) node [anchor=north west][inner sep=0.75pt]    {$a_{1}$};
\draw (272,102) node [anchor=north west][inner sep=0.75pt]    {$a_{2}$};
\draw (72,23) node [anchor=north west][inner sep=0.75pt]    {$a_{3}$};
\draw (243,102) node [anchor=north west][inner sep=0.75pt]    {$a_{3}$};
\draw (40,23) node [anchor=north west][inner sep=0.75pt]    {$a_{2}$};
\draw (102,23) node [anchor=north west][inner sep=0.75pt]    {$a_{4}$};
\draw (212,102) node [anchor=north west][inner sep=0.75pt]    {$a_{4}$};
\draw (137,23) node [anchor=north west][inner sep=0.75pt]    {$a_{k+1}$};
\draw (167,102) node [anchor=north west][inner sep=0.75pt]    {$a_{k+1}$};
\end{tikzpicture}.
\]
We compute the trace of this element, denoted by $d_a$ for short. The terms of types $2$ and $3$ vanish, because there are only pairs of elements of $H$ with intersection $0$. Thus, 

\begin{align*}
\overline{\Tr}_{2k}(d_a) &= \sum_{i=1}^{k+1}\Bigg(\begin{tikzpicture}[baseline = -16mm,x=0.75pt,y=0.75pt,yscale=-0.7,xscale=0.7]
\draw    (52.08,50) -- (52.08,110) ;
\draw    (52.08,80) -- (335.23,80) ;
\draw    (335.23,50) -- (335.23,110) ;
\draw    (100.55,50) -- (100.55,80) ;
\draw    (125.94,50) -- (125.94,80) ;
\draw    (193.65,50) -- (193.65,80) ;
\draw    (204.73,110) -- (204.73,80) ;
\draw    (245.74,110) -- (245.74,80) ;
\draw    (280.29,110) -- (280.29,80) ;
\draw [color={rgb, 255:red, 208; green, 2; blue, 27 }  ,draw opacity=1 ]   (125.94,50) .. controls (124.71,27.6) and (166,26.6) .. (166,50.6) .. controls (166,74.6) and (185.04,140.6) .. (212.12,140.6) .. controls (239.2,140.6) and (246.06,120.61) .. (245.74,110) ;
\draw  [color={rgb, 255:red, 208; green, 2; blue, 27 }  ,draw opacity=1 ] (179.88,92.97) -- (177.67,101.64) -- (171.14,95.26) ;
\draw (131.76,53.4) node [anchor=north west][inner sep=0.75pt]    {$\cdots $};
\draw (211.85,80.4) node [anchor=north west][inner sep=0.75pt]    {$\cdots $};
\draw (44.2,112) node [anchor=north west][inner sep=0.75pt]    {$a_{1}$};
\draw (327,34) node [anchor=north west][inner sep=0.75pt]    {$a_{1}$};
\draw (327,112) node [anchor=north west][inner sep=0.75pt]    {$a_{2}$};
\draw (85,34) node [anchor=north west][inner sep=0.75pt]    {$a_{i-1}$};
\draw (267.8,112) node [anchor=north west][inner sep=0.75pt]    {$a_{i-1}$};
\draw (44.2,34) node [anchor=north west][inner sep=0.75pt]    {$a_{2}$};
\draw (179.88,34) node [anchor=north west][inner sep=0.75pt]    {$a_{k+1}$};
\draw (195.58,112) node [anchor=north west][inner sep=0.75pt]    {$a_{k+1}$};
\draw (62.89,53.4) node [anchor=north west][inner sep=0.75pt]    {$\cdots $};
\draw (304.11,83.4) node [anchor=north west][inner sep=0.75pt]    {$\cdots $};
\end{tikzpicture}\Bigg)a_i
\\
&= (-1)^{k}\sum_{i=1}^{k} \overline{[...[[a_1,a_2],a_3], ..., a_{i-1}]} a_{i+1}a_{i+2}... a_{k+1}a_{k+1}a_k...a_{i+1}[...[[a_1,a_2],a_3], ..., a_{i-1}]a_i.
\end{align*}
We claim that this element is non-trivial in $C_{2k+1}(H)/\overline{C_{2k+1}(H)}$. Indeed, its reduction in $B_{2k+1}(H)$ is non-trivial. We deduce from Proposition \ref{bkspace} and its proof that the $2$-torsion of $B_{2k+1}(H)$ is a $\Z_2$-vector space of dimension $n^{k+1}$ with basis given by elements $mx\overline{m}$ for $x$ in the basis of $H$, and $m$ in the basis of $T_k(H)$ induced by the basis of $H$. In the sum above, we see that after expansion of the brackets, only one term contains $a_1$ once: $$a_2a_3...a_{k+1}a_{k+1}a_k...a_2a_1 \equiv ma_1\overline{m}$$ for $m = a_{k+1}a_k...a_2$. We can project to the line generated by this element and we get $$\overline{\Tr}(d_a) \neq 0.$$
Nevertheless, it is not hard to see that $2d_a$ is in the image of $\tau_{2k}$, by ``splitting it in half". For this, we use a vertex colored by $a_{k+2}$, hence we need $n\geq 2k+4$. We show the process for $k=3$, as the method generalizes in a straightforward way. We use the tree description of the Lie bracket of two symplectic derivations, and prove that $2d_a$ is a bracket of elements of $D_1(H)$. Since $\tau_1$ is surjective on $D_1(H)$, this proves that $2d_a \in \im(\tau_{2k})$. We have: 
$$ 2d_a = \Big[\vlfivetree{a_2}{a_1}{a_5}{a_4}{a_3} , \vlfivetree{a_4}{b_5}{a_1}{a_2}{a_3}\Big].$$ It remains to decompose these terms:

$$ \vlfivetree{a_2}{a_1}{a_5}{a_4}{a_3} = \Big[\ltritree{a_2}{a_1}{b_1} , \Big[\ltritree{a_3}{a_1}{b_1} ,\ltritree{a_4}{a_1}{a_5}\Big] \Big]$$ and we can decompose $\vlfivetree{a_4}{b_5}{a_1}{a_2}{a_3}$ similarly.
\end{proof} 

\begin{remark}
Any element in the orbit of $d_a$ under the action of the symplectic group is also a non-trivial element of torsion in the cokernel of $\tau$. From this fact we could get a (coarse) lower bound on the number of torsion elements for every even degree.
\end{remark}

\section{About a topological interpretation in degree 2, and generalizations}
\label{sec8}
Our definition of $\overline{\Tr}$ is purely algebraic. In this section, we shall give another interpretation for $\overline{\Tr}_2$, or more precisely for its reduction $\Tr_2^{\mir}$. This computation was suggested to the author by G. Massuyeau. In \cite{masY3}, Massuyeau and Meilhan introduced an invariant $\alpha$, which is the quadratic part of the relative Alexander polynomial of homology cylinders. Since it overlaps our notations, we shall use the notation $\kappa$ instead of $\alpha$.

We consider the following generalization of the Johnson filtration. For more information about the following, please refer to the survey \cite{mhsurvey}. The mapping class group can be embedded in the monoid of homology cobordisms $\mathcal{C} := \mathcal{C}_{g,1}$ over the surface $\Sigma_{g,1}$. It is known that this enlargement of $\calM$ does not act on the fundamental group of the surface, but still acts on its nilpotent quotients. The \emph{monoid of homology cylinders} $\calIC$ is the submonoid of element acting trivially on the first nilpotent quotient $H$, and we can define the Johnson filtration $\calC[k]$ as the kernel of the action on the $k$-th nilpotent quotient of $\pi$. We then generalize the definition of the Johnson homomorphisms $\tau_k$, and get the following commutative diagram:
\[
\begin{tikzcd}
\calC[k] \arrow[r, "\tau_k"]& D_{k}(H) \\
J_k \arrow[u]\arrow[ru, bend right = 10]
\end{tikzcd}.
\] An important point is that the Johnson homomorphism maps $\tau_k$ surjectively on $D_k(H)$ \cite{garlev}. We relate $\Tr^{as}$ (i.e. $\Tr_2^{\mir}$, by Proposition \ref{tracesrel}) to the invariant $\kappa : \calIC \rightarrow S_2(H)$. The map $\kappa$ is additive and is trivial on the mapping class group \cite[Prop. 3.13]{masY3}. There is a well-defined map from $S^2(H)$ to $\Lambda^2(H) \otimes \Z_2$ defined by sending $xy$ to $(x\wedge y) \otimes 1$.

\begin{proposition}
\label{propkappa}
The following diagram is commutative:
\[
\begin{tikzcd}
& \calC [2] \arrow[r, "\kappa"]\arrow[d, "\tau_2"] &S^2(H) \arrow[d] \\
J_2 \arrow[ru, bend left = 10, end anchor = {[yshift=1ex]}]\arrow[r]& D_2(H) \arrow[r, "\Tr^{as}"] &\Lambda^2(H) \otimes \Z_2.
\end{tikzcd}\]
\end{proposition}

\begin{proof}
We refer to a lot of tools defined in \cite{masY3}. According to \cite[Lemma 3.20]{masY3}, there is a commutative diagram (induced by the Birman-Craggs homomorphism of homology cylinders):  
\[
\begin{tikzcd}
\calC [3] \arrow[r]\arrow[d, "\kappa"] &H \otimes \Z_2\arrow[d] \\
S^2(H) \arrow[r] &S^2(H) \otimes \Z_2.
\end{tikzcd}\] where the map on the right is given by $x \mapsto x^2$. This precisely implies that the reduction of $\kappa$ to $\Lambda^2(H) \otimes \Z_2$ vanishes on $\calC [3]$ (because $x \wedge x = 0$). Hence, to prove the commutativity of the diagram, it is enough to prove it for generators of $\calC [2] / \calC[3] \simeq D_2(H)$.

In \cite{masY3}, Massuyeau and Meilhan give a diagrammatic interpretation of the space $\calIC[2] / Y_3$, where $Y_3$ is the so-called $Y_3$-equivalence. Equation $(4.9)$ of \cite{masY3} implies that there exists some homology cylinder $M \in \calC[2]$, such that $\tau_2(M) = \vltree{h}{i}{j}{k}$ and $$\kappa(M) \equiv \omega(h, j) ik + \omega(h, k)  ij + \omega(i, j)  hk + \omega(i, k)  hj ~~\operatorname{mod} 2.$$ This is equal to $\Tr^{as}(\vltree{h}{i}{j}{k})$ after reducing in $\Lambda^2(H) \otimes \Z_2$ (see \cite[Theorem 2.4]{faes} for the computation of $\Tr^{as}$). Furthermore, according to computations in Lemma 4.4 and to rule (G1) of the same paper, there exists a homology cylinder $M \in \calC[2]$ such that $\tau_2(M)=\frac{1}{2} \vltree{h}{k}{h}{k}$ and $$\kappa(M) = -\frac{1}{2}(\omega(h,k)2hk - 2hk) \equiv (1+\omega(h,k))hk ~~\operatorname{mod} 2.$$ Once again, this corresponds to the value of $\Tr^{as}$ after reduction.
\end{proof}

\begin{remark}
Because $\kappa$ is trivial on $J_2$, the above commutative diagram shows that $\Tr^{as}$ vanishes on $\tau_2(J_2)$. Hence we have 3 very different proofs of this fact: one using the known generators for $J_2$ (see \cite[Theorem 2.4]{faes}), one using the correspondence with $\Tr^{\mir}$, and one using the proposition above.
\end{remark}

\begin{remark}
The invariant $\kappa$ has been recently generalized by Nozaki, Sato, and Suzuki \cite{NSSred}, to a more general invariant of homology cylinders vanishing on the mapping class group. It may be interesting to investigate the relation with the Johnson homomorphisms, i.e. to generalize Proposition \ref{propkappa}, and possibly relate their invariant to $\overline{\Tr}$.
\end{remark}

\begin{remark}
\label{finalremark}
If a derivation $d$ of degree $k$ is symplectic and tame, and in the kernel of $\overline{\Tr}$, then it means by definition that there exists a lift $f_1 \in A_k$ of the derivation such that $f_1(\zeta)\zeta^{-1} \in \Gamma_{k+4}\pi$. As suggested by Heuristic C , we can continue the process by looking for a lift $f_2$ of $d$ such that $f_2(\zeta)\zeta^{-1} \in \Gamma_{k+5}\pi $ and so on. Theoretically it defines a sequence of obstructions, that might be trivial or not. Suppose all these obstructions vanish: for all $i \geq 4$, there exists $f_i \in A_k$ such that $\tau_k(f_i) = d$ and $f_i(\zeta)\zeta^{-1} \in \Gamma_{k+3+i}$. Also, we might suppose that $f_{i+1}$ was obtained as a perturbation of $f_{i}$ by multiplying by an element of $A_{k+i}$.  Unfortunately, this procedure would only yield an element in the completion $$\widehat{A}(\pi) := \lim_{\longleftarrow k} \Aut(\pi)/ A_k$$ which would be ``a lift of $d$" and ``acts trivially" on $\zeta$. By this, we mean that the completion $\widehat{A}(\pi)$ acts on all nilpotent quotients of $\pi$, and hence on the pronilpotent completion $\widehat{\pi}$ of $\pi$. Hence we can define the Andreadakis-Johnson homomorphisms on this enlargement of $\Aut(\pi)$. The elements $f_i$ induce an element $\widehat{f} \in \widehat{A}(\pi)$ such that $\tau'_k(f) = d$, and that $\widehat{f}(\zeta) = \zeta \in \widehat{\pi}$. We would like to conclude that there exists an element of $\Aut(\pi)$ with the same properties as $\widehat{f}$. Hence the question below. 
\end{remark}

\begin{question}
If we use the notation $M[k]$ for the $k$-th term of the Johnson filtration of  $M$, then the inclusion $J_k := \Aut_\zeta(\pi)[k] \subset \Aut_\zeta(\widehat{\pi})[k]$ is strict, because the Johnson homomorphisms for $\Aut_\zeta(\widehat{\pi})$ surject onto $D(H)$. Indeed, it is known that homology cylinders produce elements of $\Aut_\zeta(\widehat{\pi})$ and that the Johnson homomorphisms for homology cylinders are onto $D(H)$ in every degree \cite{garlev}. For the same reason $\widehat{A}_\zeta[k] \subset \Aut_\zeta(\widehat{\pi})[k]$ is strict. Is the inclusion $J_k \subset \widehat{A}_\zeta[k]$ strict~?
\end{question}

\begin{remark}
Because $\Tr^{as}$ is enough to compute $J_2 /J_3$, we have that $J_2/J_3 \simeq \widehat{A}_\zeta[2]/\widehat{A}_\zeta[3]$.
\end{remark}

\begin{remark}
The Satoh trace is related to the geometric constraint on the image of the Johnson homomorphisms found in \cite{KKjohnson} using the Goldman-Turaev Lie bialgebra. It would be interested to know where our new obstructions fit in this picture. Pr. Kawazumi and Pr. Sakasai suggested to the author that at least a framed version of Turaev's cobracket would be necessary to recover more than the Satoh Trace.
\end{remark}
\noindent
\bibliographystyle{plain}
\bibliography{manuscrit}
\adresse
\end{document}